\documentclass{amsart}

\newtheorem{theorem}{Theorem}[section]

\newtheorem{lemma}[theorem]{Lemma}
\newtheorem{proposition}[theorem]{Proposition}

\theoremstyle{definition}

\newtheorem{remark}[theorem]{Remark}

\newcommand{\EBF}{{\mathbf{E}}}
\newcommand{\ZBF}{{\mathbf{Z}}}
\newcommand{\NBF}{{\mathbf{N}}}
\newcommand{\RBF}{{\mathbf{R}}}
\newcommand{\CBF}{{\mathbf{C}}}
\newcommand{\eps}{{\varepsilon}}

\title[Structure, randomness, and progressions in primes]{The dichotomy between structure and randomness, arithmetic progressions, and the primes}

\author[Terence Tao]{Terence Tao}

\thanks{The author is supported by a grant from the Packard foundation.}

\begin{document}

\begin{abstract}  A famous theorem of Szemer\'edi asserts that all subsets of the integers with positive upper density will contain arbitrarily long arithmetic progressions.  There are many different proofs of this deep theorem, but they are all based on a fundamental dichotomy between structure and randomness, which in turn leads (roughly speaking) to a decomposition of any object into a structured (low-complexity) component and a random (discorrelated) component.  Important examples of these types of decompositions include the Furstenberg structure theorem and the Szemer\'edi regularity lemma.
One recent application of this dichotomy is the result of Green and Tao establishing that the prime numbers contain arbitrarily long arithmetic progressions (despite having density zero in the integers).  The power of this dichotomy is evidenced by the fact that the Green-Tao theorem requires surprisingly little technology from analytic number theory,
relying instead almost exclusively on manifestations of this dichotomy such as Szemer\'edi's theorem.  In this paper we survey various manifestations of this dichotomy in combinatorics, harmonic analysis, ergodic theory, and number theory.  As we hope to emphasize here, the underlying themes in these arguments are remarkably similar even though the contexts are radically different.
\end{abstract}



\maketitle

\section{Introduction}

In 1975, Szemer\'edi \cite{szemeredi} proved the following deep and enormously influential theorem:

\begin{theorem}[Szemer\'edi's theorem]\label{szem-thm}  Let $A$ be a subset of the integers $\ZBF$ of positive upper density, thus
$\limsup_{N \to \infty} \frac{|A \cap [-N,N]|}{|[-N,N]|} > 0$.  Here $|A|$ denotes the cardinality of a set $A$, and $[-N,N]$ denotes the integers between $-N$ and $N$.  Then for any $k \geq 3$, $A$ contains infinitely many arithmetic progressions
of length $k$.
\end{theorem}

Several proofs of this theorem are now known.  The original proof of Szemer\'edi \cite{szemeredi} was combinatorial.  A later proof of Furstenberg \cite{furst}, \cite{FKO} used ergodic theory and has led to many extensions.  A more quantitative proof of Gowers \cite{gowers-4}, \cite{gowers} was based on Fourier analysis and arithmetic combinatorics (extending a much older argument of Roth \cite{roth} handling the $k=3$ case).  A fourth proof by Gowers \cite{gowers-reg} and R\"odl, Nagle, Schacht, and Skokan \cite{nrs}, \cite{rs}, \cite{rodl}, \cite{rodl2} relied on the structural theory of hypergraphs.  These proofs are superficially all very different (with each having their own strengths and weaknesses), but have a surprising number 
of features in common.  The main difficulty in all of the proofs is that one \emph{a priori} has no control on the behaviour of the set $A$ other than a lower bound on its density; $A$ could range from being a very random set, to a very structured set, to something in between.  In each of these cases, $A$ will contain many arithmetic progressions - but the \emph{reason} for having these progressions varies from case to case.  Let us illustrate this by informally discussing some representative examples:

\begin{itemize}

\item (Random sets) Let $0 < \delta < 1$, and let $A$ be a random subset of $\ZBF$, which each integer $n$ lying in $A$ with an independent probability of $\delta$.  Then $A$ almost surely has upper density $\delta$, and it is easy to establish that $A$ almost surely has infinitely many arithmetic progressions of length $k$, basically because each progression of length $k$ in $\ZBF$ has a probability of $\delta^k$ of also lying in $A$.  A more refined version of this argument also applies when $A$ is \emph{pseudorandom} rather than random - thus we allow $A$ to be deterministic, but require that a suitable number of correlations (e.g. pair correlations, or higher order correlations) of $A$ are negligible.  The argument also extends to sparse random sets, for instance one where $\P( n \in A ) \sim 1 / \log n$.

\item (Linearly structured sets)  Consider a quasiperiodic set such as $A := \{ n: \{ \alpha n \} \leq \delta \}$,
where $0 < \delta < 1$ is fixed, $\alpha$ is a real number (e.g. $\alpha = \sqrt{2}$) and $\{x\}$ denotes the fractional part of $x$.  Such sets are ``almost periodic'' because there is a strong correlation between the events $n \in A$ and $n+L \in A$, thanks to the identity $\{\alpha (n+L) \} - \{ \alpha n \} = \{ \alpha L \}\ \mod 1$. 
An easy application of the Dirichlet approximation theorem
(to locate an approximate period $L$ with $\{ \alpha L \}$ small) shows that such sets still have infinitely many progressions of any given length $k$.  Note that this argument works regardless of whether $\alpha$ is rational or irrational.

\item (Quadratically structured sets) Consider a ``quadratically quasiperiodic'' set of the form $A := \{ n: \{ \alpha n^2 \} \leq \delta \}$.  If $\alpha$ is irrational, then this set has upper density $\delta$, thanks to Weyl's theorem on equidistribution of polynomials.  (If $\alpha$ is rational, one can still obtain some lower bound on the upper density.)  It is not linearly 
structured (there is no asymptotic correlation between the events $n \in A$ and $n+L \in A$ as $n \to \infty$ for any
fixed non-zero $L$), however it has quadratic structure in the sense that there is a strong correlation between the 
events $n \in A$, $n + L \in A$, $n + 2L \in A$, thanks to the identity
$$ \{\alpha n^2\} - 2 \{\alpha (n+L)^2\} + \{\alpha (n+2L)^2\} = 2 \{\alpha L^2\}\ \mod 1.$$
In particular $A$ does not behave like a random set.  Nevertheless, the quadratic structure still ensures that $A$ contains infinitely many arithmetic progressions of any length $k$, as one first locates a ``quadratic period'' $L$ with $\{ \alpha L^2 \}$ small, and then for suitable $n \in A$ one locates a much smaller ``linear period'' $M$ with $\{ \alpha LMn \}$ small.  If this is done correctly, the progression $n, n+LM, \ldots, n+(k-1)LM$ will be completely contained in $A$.  The same arguments also extend to a more general class of quadratically structured sets, such as the ``$2$-step nilperiodic'' set $A = \{ n: \{ \lfloor \sqrt{2} n \rfloor \sqrt{3} n \leq \delta \}$, where $\lfloor x \rfloor$ is the greatest integer function.

\item (Random subsets of structured sets) Continuing the previous example $A := \{ n: \{ \alpha n^2 \} \leq \delta \}$, let
$A'$ be a random subset of $A$ with each $n \in A$ lying in $A'$ with an independent probability of $\delta'$ for some $0 < \delta' < 1$.  Then this set $A'$ almost surely has a positive density of $\delta \delta'$ if $\alpha$ is irrational.  The set $A'$ almost surely has infinitely many progressions of length $k$, since $A$ already starts with infinitely many such progressions, and each such progression as a probability of $(\delta')^k$ of also lying in $A'$.  One can generalize this example to random sets $\tilde A$ where the events $n \in \tilde A$ are independent as $n$ varies, and the probability $P(n \in \tilde A)$ is a ``quadratically almost periodic'' function of $n$ such as $P(n \in \tilde A) = F( \{ \alpha n^2 \} )$ for some nice (e.g. piecewise continuous) function $F$ taking values between $0$ and $1$; the preceding example is the case where $F(x) := \delta' 1_{x < \delta}$. It is also possible to adapt this argument to (possibly sparse) pseudorandom subsets of structured sets, though one needs to take some care in defining exactly what ``pseudorandom'' means here.  

\item (Sets containing random subsets of structured sets) Let $A''$ be any set which contains the set $A'$ (or $\tilde A$) of the previous example.  Since $A'$ contains infinitely many progressions of length $k$, it is trivial that $A''$ does also.

\end{itemize}

As the above examples should make clear, the reason for the truth of Szemer\'edi's theorem is very different in the cases when $A$ is random, and when $A$ is structured.  These two cases can then be combined to handle the case when $A$ is (or contains) a large (pseudo-)random subset of a structured set.  Each of the proofs of Szemer\'edi's theorem now hinge on a \emph{structure 
theorem} which, very roughly speaking, asserts that \emph{every} set of positive density is (or contains) a large pseudorandom subset of a structured set; each of the four proofs obtains a structure theorem of this sort in a different way (and in a very different language).  These remarkable structural results - which include the Furstenberg structure theorem and the Szemer\'edi regularity lemma as examples - are of 
independent interest (beyond their immediate applications to arithmetic progressions), and have led to many further developments and insights.  For instance, in \cite{gt-primes} a ``weighted'' structure theorem
(which was in some sense a hybrid of the Furstenberg structure theorem and the Szemer\'edi regularity lemma) was the primary new ingredient in proving that the primes $P := \{2,3,5,7,\ldots\}$ contained arbitrarily long arithmetic progressions.  While that latter claim is ostensibly a number-theoretical result, the method of proof in fact uses surprisingly little from number theory, being much closer in spirit to the proofs of Szemer\'edi's theorem (and in fact Szemer\'edi's theorem is a crucial ingredient in the proof).  This can be seen from the fact that the argument in \cite{gt-primes} in fact proves the following stronger result:

\begin{theorem}[Szemer\'edi's theorem in the primes]\label{szemp-thm}\cite{gt-primes}  Let $A$ be a subset of the primes $P$ of positive \emph{relative} upper density, thus
$\limsup_{N \to \infty} \frac{|A \cap [-N,N]|}{|P \cap [-N,N]|} > 0$.  Then for any $k \geq 3$, $A$ contains 
infinitely many arithmetic progressions of length $k$.
\end{theorem}

This result was first established in the $k=3$ case by Green \cite{green}, the key step again being a (Fourier-analytic) 
structure theorem, this time for subsets of the primes.  The arguments used to prove this theorem do not directly address 
the important question of whether the primes $P$ (or any subset thereof) have any pseudorandomness properties (but see Section \ref{primes-sec} below).  However, the structure theorem does allow one to (essentially) describe any dense subset of the primes as
a (sparse) pseudorandom subset of some unspecified dense set, which turns out to be sufficient (thanks to Szemer\'edi's theorem) for the purpose of establishing the existence of arithmetic progressions.

There are now several expositions of Theorem \ref{szemp-thm}; see for instance \cite{kra-survey}, \cite{green-survey}, \cite{tao-coates}, \cite{tao-survey}, \cite{host-survey}.  Rather than give another exposition of this result, we have chosen to take a broader view, surveying the collection of structural theorems which underlie the proof of such results as 
Theorem \ref{szem-thm} and Theorem \ref{szemp-thm}.  These theorems have remarkably varied contexts - measure theory, ergodic theory, graph theory, hypergraph theory, probability theory, information theory, and Fourier analysis - and can be either qualitative (infinitary) or quantitative (finitary) in nature.  However, their \emph{proofs} tend to share a number of common features, and thus serve as a kind of ``Rosetta stone'' connecting these various fields.  Firstly, for a given class of objects, one quantifies what it means for an object to be ``(pseudo-)random'' and an object to be ``structured''.  Then, one establishes a \emph{dichotomy between randomness and structure}, which typically looks something like this:

\begin{quotation}
\emph{If an object is not (pseudo-)random, then it (or some non-trivial component of it) correlates with a structured object.}
\end{quotation}

One can then iterate this dichotomy repeatedly (e.g. via a stopping time argument, or by Zorn's lemma), to extract out all the correlations with structured objects, to obtain a \emph{weak structure theorem} which typically looks as follows:

\begin{quotation}
\emph{If $A$ is an arbitrary object, then $A$ (or some non-trivial component of $A$) splits as the sum of a structured object, plus a pseudorandom error.}
\end{quotation}

In many circumstances, we need to improve this result to a \emph{strong structure theorem}:

\begin{quotation}
\emph{If $A$ is an arbitrary object, then $A$ (or some non-trivial component of $A$) splits as the sum of a structured object, plus a small error, plus a \textbf{very} pseudorandom error.}
\end{quotation}

When one is working in an infinitary (qualitative) setting rather than a finitary (quantitative) one - which is for instance the case in the ergodic theory approach - one works instead with an \emph{asymptotic structure theorem}:

\begin{quotation}
\emph{If $A$ is an arbitrary object, then $A$ (or some non-trivial component of $A$) splits as the sum of a ``compact'' object (the limit of structured objects), plus an \textbf{infinitely} pseudorandom error.}
\end{quotation}

The reason for the terminology ``compact'' to describe the limit of structured objects is in analogy to how a compact operator can be viewed as the limit of finite rank operators; see \cite{furst-book} for further discussion.

In many applications, the small or pseudorandom errors in these structure theorems are negligible, and one then reduces to the study of structured objects.  One then exploits the structure of these objects to conclude the desired application.

Our focus here is on the structure theorems related to Szemer\'edi's theorem and related results such as 
Theorem \ref{szemp-thm}; we will not have space to describe all the generalizations and refinements of these results here.  However, these types of structural theorems appear in other contexts also, for instance the Koml\'os subsequence principle \cite{komlos-subseq} in 
probability theory.  The Lebesgue decomposition of a spectral measure into pure point, singular continuous, and absolutely continuous spectral components can also be viewed as a structure theorem of the above type.  Also, the stopping time arguments which underlie the structural theorems here are also widely used in harmonic analysis, in particular obtaining fundamental decompositions such as the Calder\'on-Zygmund decomposition or the atomic decomposition of Hardy spaces (see e.g. \cite{stein:large}), as well as the tree selection arguments used in multilinear harmonic analysis (see e.g. \cite{laceyt1}).
It may be worth investigating whether there are any concrete connections between these disparate structural theorems.

\section{Ergodic theory}

We now illustrate the above general strategy in a number of contexts, beginning with the ergodic theory approach to Szemer\'edi's
theorem, where the dichotomy between structure and randomness is particularly clean and explicit, and one can work with an asymptotic structure theorem rather than a weak or strong one.  Very informally speaking,
the ergodic theory approach seeks to understand the set $A$ of integers by analyzing the asymptotic correlations of the shifts
$A+n := \{ a+n: a \in A \}$ (or of various asymptotic averages of these shifts), and treating these shifts as occurring on an
abstract measure space.
More formally, let $X$ be a measure space with probability measure $d\mu$, and 
let $T: X \to X$ be a bijection such that $T$ and $T^{-1}$ are both measure-preserving maps.  The associated shift operator
$T: f \mapsto f \circ T^{-1}$ is thus a unitary operator on the Hilbert space $L^2(X)$ of complex-valued square-integrable functions with the usual inner product $\langle f, g\rangle :=\int_X f\overline{g}\ d\mu$.  A famous transference result known as 
the \emph{Furstenberg correspondence principle}\footnote{Morally speaking, to deduce Szemer\'edi's theorem from Furstenberg's theorem, one takes $X$ to be the integers $\ZBF$, $T$ to be the standard shift $n \mapsto n+1$, and $\mu$ to be the density
$\mu(A) = \lim_{N \to \infty} \frac{|A \cap [-N,N]|}{|[-N,N]|}$.  This does not quite work because not all sets $A$ have a well-defined density, however additional arguments (e.g. using the Hahn-Banach theorem) can fix this problem.}
 (see \cite{furst}, \cite{FKO}, \cite{furst-book}) shows that Szemer\'edi's theorem is then equivalent to

\begin{theorem}[Furstenberg recurrence theorem]\label{furst-thm}\cite{furst}  Let $X$ and $T$ be as above, and let $f \in L^\infty(X)$ be any bounded non-negative function with $\int_X f\ d\mu > 0$.  Then for any $k \geq 1$ we have
$$ \liminf_{N \to \infty} \EBF_{1 \leq n \leq N} \int_X f T^n f \ldots T^{(k-1)n} f\ d\mu > 0.$$
Here and in the sequel we use $\EBF_{n \in I} a_n$ as a shorthand for the average $\frac{1}{|I|}\sum_{n \in I} a_n$.
\end{theorem}

When $k=2$ this is essentially the Poincar\'e recurrence theorem; by using the von Neumann ergodic theorem one can also show that the limit exists (thus the lim inf can be replaced with a lim).  
The $k=3$ case can be proved by the following argument, as observed in \cite{furst-book}. We need to show that
\begin{equation}\label{triple}
\liminf_{N\to \infty} \EBF_{1 \leq n \leq N} \int_X f T^n f T^{2n} f\ d\mu > 0
\end{equation}
whenever $f$ is bounded, non-negative, and has positive integral.

The first key observation is that any sufficiently pseudorandom component of $f$ will give a negligible contribution 
to \eqref{triple} and can be dropped.  More precisely, let us call  $f$ is \emph{linearly pseudorandom} 
(or \emph{weakly mixing}) with respect to the shift $T$ if
we have 
\begin{equation}\label{unif-def}
\lim_{N \to\infty} \EBF_{1 \leq n \leq N} |\langle T^n f, f \rangle|^2 = 0.
\end{equation}
Such functions are negligible for the
purpose of computing averages such as those in \eqref{triple}; indeed, if at least one of $f, g, h \in L^\infty(X)$ is 
linearly pseudorandom, then an easy application of van der Corput's lemma (which in turn is an application of Cauchy-Schwarz) shows that 
$$\lim_{N \to \infty}\EBF_{1 \leq n \leq N} \int_X f T^n g T^{2n} h\ d\mu = 0.$$
We shall refer to these types of results - that pseudorandom functions are negligible when averaged against other functions - as \emph{generalized von Neumann theorems}.

In view of this generalized von Neumann theorem, one is now tempted to ``quotient out'' all the pseudorandom functions
and work with a reduced class of ``structured'' functions.  In this particular case, it turns out that the correct notion of structure is that of a \emph{linearly almost periodic function}, which are in turn generated by the \emph{linear 
eigenfunctions} of $T$.  To
make this more precise, we need the following dichotomy:

\begin{lemma}[Dichotomy between randomness and structure]\label{lemcor} Suppose that $f \in L^\infty(X)$ is not linearly pseudorandom.  Then there exists an linear eigenfunction $g \in L^\infty(X)$ of $T$ (thus $Tg=\lambda g$ for some $\lambda \in \CBF$) such that $\langle f, g \rangle \neq 0$.
\end{lemma}

\begin{remark} Observe that if $g$ is a linear eigenfunction of $T$ with $Tg = \lambda g$, then $|\lambda|=1$ and $\lim_{N \to \infty} \EBF_{1\leq n \leq N} \int_X g T^n \overline{g}^2 T^{2n} g\ d\mu = \int_X |g|^4$.  Thus linear eigenfunctions can and do give nontrivial contributions to the expression in \eqref{triple}. One can view Lemma \ref{lemcor} as a converse to this observation.
\end{remark}

\begin{proof}(Sketch) Let $S$ denote the operator $S g := \lim_{N \to \infty} \EBF_{1\leq n \leq N} \langle g, T^n f \rangle T^n f$ (this limit exists by the von Neumann ergodic theorem).  One can show that $S$ is self-adjoint, compact, and commutes with $T$, and thus by spectral theory has an expansion of the form $S g = \sum_k c_k \langle g, g_k \rangle g_k$ where $g_k$ are a countable sequence of eigenfunctions of $T$ and $c_k$ are scalars.  Since $f$ is not linearly pseudorandom, we have $\langle Sf, f \rangle > 0$, so in particular $Sf$ is non-zero.  This implies that $\langle f, g_k \rangle \neq 0$ for one of the eigenfunctions $g_k$, and we are done.  (The eigenfunctions must be bounded since $S$ maps $L^2(X)$ to $L^\infty(X)$.)
\end{proof}

This lemma has the following consequence.  Let ${\mathcal Z}_1$ be the $\sigma$-algebra generated by all the eigenfunctions of $T$, this is known as the \emph{Kronecker factor} of $X$, and roughly speaking encapsulates all the ``linear structure'' in the measure preserving system.  Given every function $f \in L^2(X)$, we have the decomposition $f = f_{U^\perp} + f_U$, where $f_{U^\perp} := \EBF(f|{\mathcal Z}_1)$ is the conditional expectation of $f$ with respect to the $\sigma$-algebra ${\mathcal Z}_1$ (i.e. the orthogonal projection from $L^2(X)$ to the ${\mathcal Z}_1$-measurable functions). By construction, $f_U := f - \EBF(f|{\mathcal Z}_1)$ is orthogonal to every eigenfunction of $T$, and is hence linearly pseudorandom by Lemma \ref{lemcor}.  In particular, we have established 

\begin{proposition}[Asymptotic structure theorem]\label{struct-3} Let $f$ be bounded and non-negative, with positive integral.  Then we can split\footnote{The notation is from \cite{gt-primes}; the subscript $U$ stands for ``Gowers uniform'' (pseudorandom), and $U^\perp$ for ``Gowers anti-uniform'' (structured).} $f = f_{U^\perp} + f_U$, where $f_{U^\perp}$ is bounded, non-negative, and ${\mathcal Z}_1$-measurable (and thus approximable in $L^2$ to arbitrary accuracy by finite linear combinations of linear eigenfunctions), with positive integral, and $f_U$ is linearly pseudorandom.
\end{proposition}

This result is closely related to the Koopman-von Neumann theorem in ergodic theory.  In the language of the introduction, it asserts (very roughly speaking) that any set $A$ of integers can be viewed as a (linearly) pseudorandom set where the ``probability'' $f_{U^\perp}(n)$ that a given element $n$ lies in $A$ is a (linearly) almost periodic function of $n$. 

Note that the linearly pseudorandom component $f_U$ of $f$ gives no contribution to \eqref{triple}, thanks to the generalized von Neumann theorem.  Thus we may freely replace $f$ by $f_{U^\perp}$ if desired; in other words, for the purposes of proving \eqref{triple} we may assume without loss of generality that $f$ is measurable with respect to the Kronecker factor ${\mathcal Z}_1$.  In the notation of \cite{furst-weiss}, we have just shown that the Kronecker factor is a \emph{characteristic factor} for 
the recurrence in \eqref{triple}.  (In fact it is essentially the universal factor for this recurrence, see \cite{ziegler}, \cite{host-kra2} for further discussion.)

We have reduced the proof of \eqref{triple} to the case when $f$ is structured, in the sense of being measurable 
in ${\mathcal Z}_2$.  There are two ways to obtain the desired ``structured recurrence'' result.  Firstly there is a ``soft'' approach, in which one observes that every ${\mathcal Z}_1$-measurable square-integrable function $f$ is \emph{almost periodic}, in the sense that for any $\eps > 0$ there exists a set of integers $n$ of positive density such that $T^n f$ is within $\eps$ of $f$ in $L^2(X)$; from this it is easy to show that $\int_X f T^n f T^{2n} f\ d\mu$ is close to $\int_X f^3$ for a set of integers $n$ of positive density, which implies \eqref{triple}.  This almost periodicity can be verified by first checking it for polynomial combinations of linear eigenfunctions, and then extending by density arguments.  There is also a ``hard'' approach, in which one obtains algebraic and topological control on the Kronecker factor ${\mathcal Z}_1$.  In fact, from a spectral analysis of $T$ one can show that ${\mathcal Z}_1$ is the inverse limit of a sequence of $\sigma$-algebras, on each of which the shift $T$ is isomorphic to a 
shift $x \mapsto x+\alpha$ on a compact abelian Lie group $G$.  This gives a very concrete description of the functions $f$ which are measurable in the Kronecker factor, and one can establish \eqref{triple} by a direct argument similar to that used in
in the introduction for linearly structured sets.  This ``hard'' approach gives a bit more information; for instance, it can be used to show that the limit in \eqref{triple} actually converges, so one can replace the lim inf by a lim.

It turns out that these arguments extend (with some non-trivial effort) to the case of higher $k$.  For sake of exposition let us just discuss the $k=4$ case, though most of the assertions here extend to higher $k$.  We wish to prove that
\begin{equation}\label{quadruple}
\liminf_{N\to \infty} \EBF_{1 \leq n \leq N} \int_X f T^n f T^{2n} f T^{3n} f\ d\mu > 0
\end{equation}
whenever $f$ is bounded, non-negative, and has positive integral.  Here, it turns out that we must strengthen the notion of pseudorandomness (and hence generalize the notion of structure); linear pseudorandomness is no longer sufficient to imply negligibility.  For instance, let $f$ be a \emph{quadratic eigenfunction}, in the sense that $Tf = \lambda f$, where $\lambda$
is no longer constant but is itself a linear eigenfunction, thus $T\lambda = c\lambda$ for some constant $c$. As an example,
if $X = (\RBF/\ZBF)^2$ with the skew shift $T(x,y) = (x+\alpha,y+x)$ for some fixed number $\alpha$, then the function $f(x,y) = e^{2\pi i y}$ is a quadratic eigenfunction but not a linear one.  Typically such quadratic eigenfunctions will be linearly pseudorandom, but if $|\lambda|=|c|=1$ (which is often the case) then we have the identity
\begin{equation}\label{erg-ex}
 \EBF_{1 \leq n \leq N} \int_X f T^n \overline{f}^3 T^{2n} f^3 T^{3n} \overline{f}\ d\mu = \int_X |f|^8\ d\mu
 \end{equation}
and so we see that these functions can give non-trivial contributions to expressions such as \eqref{triple}.  The correct notion of pseudorandomness is now \emph{quadratic pseudorandomness}, by which we mean that
$$ \lim_{H \to \infty} \lim_{N \to \infty} \EBF_{1 \leq n\leq N} \EBF_{1 \leq h \leq H} |\langle T^h f \overline{f}, T^n( T^h f \overline{f} ) \rangle|^2 = 0.$$
In other words, $f$ is quadratically pseudorandom if and only if 
$T^h f \overline{f}$ is asymptotically linearly pseudorandom on the average
as $h \to \infty$.  Several applications of van der Corput's lemma give a generalized von Neumann theorem, asserting that
$$ \lim_{N\to \infty} \EBF_{1 \leq n \leq N} \int_X f_0 T^n f_1 T^{2n} f_2 T^{3n} f_3\ d\mu = 0$$
whenever $f_0,f_1,f_2,f_3$ are bounded functions with at least one function quadratically pseudorandom.

One would now like to construct a factor ${\mathcal Z}_2$ (presumably larger than the Kronecker factor ${\mathcal Z}_1$) which will play the role of the Kronecker factor for the average \eqref{quadruple}; in particular, we would like a statement of the form

\begin{lemma}[Dichotomy between randomness and structure]\label{lemcor-4} Suppose that $f \in L^\infty(X)$ is not linearly pseudorandom.  Then there exists a ${\mathcal Z}_2$-measurable function $g \in L^\infty(X)$ such that $\langle f, g \rangle \neq 0$.
\end{lemma}

which would imply\footnote{One can generalize this structure theorem to obtain similar characteristic factors ${\mathcal Z}_3$, ${\mathcal Z}_4$ for cubic pseudorandomness, quartic pseudorandomness, etc.  Applying Zorn's lemma, one eventually obtains the \emph{Furstenberg structure theorem}, which decomposes any measure preserving system as a weakly mixing extension of a distal system, and thus decomposes any function as a distal function plus an ``infinitely pseudorandom'' error; see \cite{FKO}.  However this decomposition is not the most ``efficient'' way to prove Szemer\'edi's theorem, as the notion of pseudorandomness is too strong, and hence the notion of structure too general.  It does illustrate however that one does have considerable flexibility in where to draw the line between randomness and structure.}

\begin{proposition}[Asymptotic structure theorem]\label{struct-4} Let $f$ be bounded and non-negative, with positive integral.  Then we can split $f = f_{U^\perp} + f_U$, where $f_{U^\perp}$ is bounded, non-negative, and ${\mathcal Z}_2$-measurable, with positive integral,and $f_U$ is quadratically pseudorandom.
\end{proposition}

This reduces the proof of \eqref{quadruple} to that of ${\mathcal Z}_2$-measurable $f$. The existence of such a factor ${\mathcal Z}_2$ (which would be a \emph{characteristic factor} for this average) is trivial to construct, as we could just take ${\mathcal Z}_2$ to be the entire $\sigma$-algebra, and it is in fact easy (via Zorn's lemma) to show the existence of a ``best'' such factor, which embed into all other characteristic factors for this average (see \cite{ziegler}). Of course, for the concept of characteristic factor to be useful we would like ${\mathcal Z}_2$ to be as small as possible, and furthermore to have some concrete structural description of the factor.  An obvious guess
for ${\mathcal Z}_2$ would be the $\sigma$-algebra generated by all the linear and quadratic eigenfunctions, but this factor turns out to be a bit too small (see \cite{furst-weiss}; this is related to the example of the $2$-step nilperiodic set in the introduction).  A more effective candidate for ${\mathcal Z}_2$, analogous to the ``soft'' description of the Kronecker factor, is the space of all ``quadratically almost periodic functions''.  This concept is a bit tricky to define rigorously (see e.g. \cite{FKO}, \cite{furst-book}, \cite{tao:szemeredi}), but roughly speaking, a function $f$ is linearly almost periodic if the orbit
$\{ T^n f: n \in \ZBF\}$ is precompact in $L^2(X)$ viewed as a Hilbert space, while a function $f$ is quadratically almost periodic if the orbit is precompact in $L^2(X)$ viewed as a Hilbert \emph{module} over the Kronecker factor $L^\infty({\mathcal Z}_1)$; this can be viewed as a matrix-valued (or more precisely compact operator-valued) extension of the concept of a quadratic eigenfunction.
Another rough definition is as follows: a function $f$ is linearly almost periodic if $T^n f(x)$ is close to $f(x)$ for many constants $n$, whereas a function $f$ is quadratically almost periodic if $T^{n(x)} f(x)$ is close to $f(x)$ for a function $n(x)$ which is itself linearly almost periodic.  It turns out that with this ``soft'' proposal for ${\mathcal Z}_2$,
it is easy to prove Lemma \ref{lemcor-4} and hence Proposition \ref{struct-4}, essentially by obtaining a ``relative'' version of
the proof of Lemma \ref{lemcor}.  The derivation of \eqref{quadruple} in this soft factor is slightly tricky though, requiring either van der Waerden's theorem, or the color focusing argument used to prove van der Waerden's theorem; see \cite{furst}, \cite{FKO}, \cite{furst-book}, \cite{tao:szemeredi}.  
More recently, a more efficient ``hard'' factor ${\mathcal Z}_2$ was constructed by Conze-Lesigne \cite{cl}, Furstenberg-Weiss \cite{furst-weiss}, and Host-Kra \cite{hk1}; the analogous factors for higher $k$ are more difficult to construct, but this was achieved by Host-Kra in \cite{host-kra2}, and also subsequently by Ziegler \cite{ziegler}.  This factor yields more precise information, including convergence of the limit in \eqref{quadruple}.  Here, the concept of a \emph{$2$-step nilsystem} is used to define structure.  A $2$-step nilsystem is a compact symmetric space $G/\Gamma$, with $G$ a $2$-step nilpotent Lie group and $\Gamma$ is a closed subgroup, together with a shift element $\alpha \in G$, which generates a shift $T(x\Gamma) := \alpha x\Gamma$.  The factor ${\mathcal Z}_2$ constructed in these papers is then the inverse limit of a sequence of $\sigma$-algebras, on which the shift is equivalent to a $2$-step nilsystem.  This should be compared with the ``hard'' description of the Kronecker factor, which is the $1$-step analogue of the above result.  Establishing the bound \eqref{quadruple} then reduces to the problem of understanding the structure of arithmetic progressions $x\Gamma$, $\alpha x \Gamma$, $\alpha^2 x \Gamma$, $\alpha^3 x \Gamma$ on the nilsystem, which can be handled by algebraic arguments, for instance using the machinery of Hall-Petresco sequences \cite{leibman}.

The ergodic methods, while non-elementary and non-quantitative (though see \cite{tao:szemeredi}), have proven to be the 
most powerful and flexible approach to Szemer\'edi's theorem, leading to many generalizations and refinements.  However, it seems
that a purely ``soft'' ergodic approach is not quite capable by itself of extending to the primes as in Theorem \ref{szemp-thm}, though it comes tantalizingly close.  In particular, one can use Theorem \ref{furst-thm} and a variant of the Furstenberg correspondence principle to establish Theorem \ref{szemp-thm} when the set of primes $P$ is replaced by a random subset $\tilde P$ of the positive integers, with $n \in \tilde P$ with independent probability $1/\log n$ for $n > 1$; see \cite{tao:correspond}.
Roughly speaking, if $A$ is a subset of $\tilde P$, the idea is to construct an abstract measure-preserving system generated by a set $\tilde A$, in which $\mu(T^{n_1} \tilde A \cap\ldots \cap T^{n_k} \tilde A)$ is the normalized density of
$(A+n_1) \cap \ldots \cap (A+n_k)$ for any $n_1,\ldots,n_k$.  Unfortunately, this approach requires the ambient space $\tilde P$ to be extremely pseudorandom and does not seem to extend easily to the primes.

\section{Fourier analysis}

We now turn to a more quantitative approach to Szemer\'edi's theorem, based primarily on Fourier analysis and arithmetic combinatorics.  Here, one analyzes a set of integers $A$ finitarily, truncating to a finite setting such as the discrete integral $\{1,\ldots,N\}$ or the cyclic group $\ZBF/N\ZBF$, and then testing the correlations of $A$ with linear phases such as $n \mapsto e^{2\pi i k n / N}$, quadratic phases $n \mapsto e^{2\pi i k n^2/N}$, or similar objects.  This approach has lead to the best known bounds on Szemer\'edi's theorem, though it has not yet been able to handle many of the generalizations of this theorem that can
be treated by ergodic or graph-theoretic methods.  In analogy with the ergodic arguments,
the $k=3$ case of Szemer\'edi's theorem can be handled by linear Fourier analysis (as was done by Roth \cite{roth}), while the $k=4$ case requires quadratic Fourier analysis (as was done by Gowers \cite{gowers-4}), and so forth for higher order $k$ (see \cite{gowers}).  The Fourier analytic approach seems to be closely related to the theory of the ``hard'' characteristic factors discovered in the ergodic theory arguments, although the precise nature of this relationship is still being understood.

It is convenient to work in a cyclic group $\ZBF/N\ZBF$ of prime order.  It can be shown via averaging arguments (see \cite{varnavides}) that Szemer\'edi's theorem is equivalent to the following quantitative version:

\begin{theorem}[Szemer\'edi's theorem, quantitative version]\label{sz-quant} Let $N > 1$ be a large prime, let $k \geq 3$, and let $0 < \delta < 1$.  Let $f: \ZBF/N\ZBF \to \RBF$ be a function with $0 \leq f(x)\leq 1$ for all $x \in \ZBF/N\ZBF$ and $\EBF_{x \in \ZBF/N\ZBF} f(x) \geq \delta$.  Then we have
$$ \EBF_{x,r \in \ZBF/N\ZBF} f(x) T^r f(x) \ldots T^{(k-1)r} f(x) \geq c(k,\delta)$$
for some $c(k,\delta) > 0$ depending only on $k$ and $\delta$, where $T^r f(x) := f(x+r)$ is the shift operator on $\ZBF/N\ZBF$.
\end{theorem}

We remark that the Fourier-analytic arguments in Gowers \cite{gowers} give the best known lower bounds on $c(k,\delta)$, namely
$c(k,\delta) > 2^{-2^{1/\delta^{c_k}}}$ where $c_k := 2^{2^{k+9}}$.  In the $k=3$ case it is known that $c(3,\delta) \geq \delta^{C/\delta^2}$ for some absolute constant $C$, see \cite{bourgain-triples}.  A conjecture of Erd\H{o}s and Tur\'an \cite{erdos} is roughly equivalent to asserting that $c(k,\delta) > e^{-C_k/\delta}$ for some $C_k$.  In the converse direction,
an example of Behrend shows that $c(3,\delta)$ cannot exceed $e^{c \log^2(1/\delta)}$ for some small absolute constant $c$, with similar results for higher values of $k$; in particular, $c(k,\delta)$ cannot be as large as any fixed power of $\delta$.  This already rules out a number of elementary approaches to Szemer\'edi's theorem and suggests that any proof must involve some sort of iterative argument.  

Let us first describe (in more ``modern'' language) Roth's original proof \cite{roth} of Szemer\'edi's theorem in the $k=3$ case.  We need to establish a bound of the form
\begin{equation}\label{triple-fourier}
\EBF_{x,r \in \ZBF/N\ZBF} f(x) T^r f(x) T^{2r} f(x) \geq c(3,\delta) > 0
\end{equation}
when $f$ takes values between $0$ and $1$ and has mean at least $\delta$.  As in the ergodic argument, we first look for a notion of pseudorandomness which will ensure that the average in \eqref{triple-fourier} is negligible.  It is convenient to introduce the \emph{Gowers $U^2(\ZBF/N\ZBF)$ uniformity norm} by the formula
$$ \| f \|_{U^2(\ZBF/N\ZBF)}^4 := \EBF_{n \in \ZBF/N\ZBF} |\EBF_{x \in \ZBF/N\ZBF} T^n f(x) \overline{f(x)}|^2,$$
and informally refer to $f$ as \emph{linearly pseudorandom} (or \emph{linearly Gowers-uniform}) if its $U^2$ norm is small; compare this with \eqref{unif-def}.  The $U^2$ norm is indeed a norm; this can be verified either by several applications of the Cauchy-Schwarz inequality, or via the Fourier identity 
\begin{equation}\label{fourier-form}
\| f \|_{U^2(\ZBF/N\ZBF)}^4 = \sum_{\xi \in \ZBF/N\ZBF} |\hat f(\xi)|^4,
\end{equation}
where
$\hat f(\xi) := \EBF_{x \in \ZBF/N\ZBF} f(x) e^{-2\pi i x\xi / N}$ is the usual Fourier transform.  Some further applications of Cauchy-Schwarz (or Plancherel's theorem and H\"older's inequality) yields the generalized von Neumann theorem
\begin{equation}\label{gvn}
|\EBF_{x,r \in \ZBF/N\ZBF} f_0(x) T^r f_1(x) T^{2r} f_2(x)| \leq \min_{j=0,1,2}  \|f_j\|_{U^2(\ZBF/N\ZBF)}
\end{equation}
whenever $f_0,f_1,f_2$ are bounded in magnitude by $1$.  Thus, as before, linearly pseudorandom functions give a small contribution to the average in \eqref{triple-fourier}, though now that we are in a finitary setting the contribution does not vanish completely.

The next step is to establish a dichotomy between linear pseudorandomness and some sort of usable structure.  From \eqref{fourier-form} and Plancherel's theorem we easily obtain the following analogue of Lemma \ref{lemcor}:

\begin{lemma}[Dichotomy between randomness and structure]\label{lemcor-planch} Suppose that $f: \ZBF/N\ZBF \to \CBF$ is bounded in magnitude by $1$ with $\|f\|_{U^2(\ZBF/N\ZBF)} \geq \eta$ for some $0 < \eta < 1$.  Then there exists a linear phase function $\phi: \ZBF/N\ZBF \to \RBF/\ZBF$ (thus $\phi(x) = \xi x/N + c$ for some $\xi \in \ZBF/N\ZBF$ and $c \in \RBF/\ZBF$) such that
$|\EBF_{x \in \ZBF/N\ZBF} f(x) e^{-2\pi i \phi(x)}| \geq \eta^2$.
\end{lemma}

The next step is to iterate this lemma to obtain a suitable structure theorem.  There are two slightly different ways to do this.
Firstly there is the original \emph{density increment argument} approach of Roth \cite{roth}, which we sketch as follows.
It is convenient to work on a discrete interval $[1,N/3]$, which we identify with a subset of $\ZBF/N\ZBF$ in the obvious manner.  Let $f: [1,N/3] \to \RBF$ be a non-negative function bounded in magnitude by $1$, and let $\eta$ be a parameter to be
chosen later.  If $f-\EBF_{1 \leq x \leq N/3} f(x)$ is not linearly pseudorandom, in the sense that $\|f - \EBF_{1 \leq x \leq N/3} f(x)\|_{U^2(\ZBF/N\ZBF)} \geq \eta$, then we apply
Lemma \ref{lemcor-planch} to obtain a correlation with a linear phase $\phi$.  An easy application of the Dirichlet approximation theorem then shows that one can partition $[1,N/3]$ into arithmetic progressions (of length roughly $\eta^2 \sqrt{N}$) on which $\phi$ is essentially 
constant (fluctuating by at most $\eta^2/100$, say).  A pigeonhole argument (exploiting the fact that $f - \EBF_{1 \leq x \leq N/3} f(x)$ has mean zero) then shows that on one of these progressions, say $P$,
$f$ has significantly higher density than on the average, in the sense that $\EBF_{x \in P} f(x) \geq \EBF_{x \in \ZBF/N\ZBF} f(x) + \eta^2/100$.  One can then apply an affine transformation to convert this progression $P$ into another discrete interval $\{1,\ldots,N'/3\}$, where $N'$ is essentially the square root of $N$.  One then iterates this argument until linear pseudorandomness is obtained (using the fact that the density of $f$ cannot increase beyond $1$), and one eventually obtains

\begin{theorem}[Structure theorem]\label{stru}  Let $f: [1,N/3] \to \RBF$ be a non-negative function bounded by $1$, and let $\eta > 0$.  Then there exists a progression $P$ in $[1,N/3]$ of length at least $c(\eta) N^{c(\eta)}$ for some $c(\eta) > 0$, on which we have the splitting $f = f_{U^\perp} + f_U$, where $f_U^{\perp} := \EBF_{x \in P} f(x) \geq \EBF_{1 \leq x \leq N/3} f(x)$
is the mean of $f$ on $P$, and $f_U$ is linearly pseudorandom in the sense that
$$ \| f_U \|_{U^2(\ZBF/M\ZBF)} \leq \eta$$
where we identify $P$ with a subset of a cyclic group $\ZBF/M\ZBF$ of cardinality $M \approx 3|P|$ in the usual manner.
\end{theorem}

More informally, any function will contain an arithmetic  progression $P$ of significant size on which $f$ can be decomposed into
a non-trivial structured component $f_{U^\perp}$ and a pseudorandom component $f_U$.  In the language of the introduction, it is essentially saying that any dense set $A$ of integers will contain components which are dense pseudorandom subsets of long progressions.  Once one has this theorem, it is an easy matter to establish Szemer\'edi's theorem in the $k=3$ case.  Indeed, if $A \subseteq \ZBF$ has upper density greater than $\delta$, then we can find arbitrarily large primes $N$ such that $|A \cap [1,N/3]|
\geq \delta N/3$.  Applying Theorem \ref{stru} with $\eta := \delta^3/100$, and $f$ equal to the indicator function of $A \cap [1,N/3]$, we can find a progression $P$ in $\{1,\ldots,N/3\}$ of length at least $c(\delta) N^{c(\delta)}$ on which $\EBF_{x \in P} f(x) \geq \delta$ and $f - \EBF_{x \in P} f(x)$ is linearly pseudorandom in the sense of Theorem \ref{stru}.  It is then an easy matter to apply the generalized von Neumann theorem to show that $A \cap P$ contains many arithmetic progressions of length
three (in fact it contains $\gg \delta^3 |P|^3$ such progressions).  Letting $N$ (and hence $|P|$) tend to infinity we obtain
Szemer\'edi's theorem in the $k=3$ case.  An averaging argument of Varnavides \cite{varnavides} then yields the more quantitative version in Theorem \ref{sz-quant} (but with a moderately bad bound for $c(3,\delta)$, namely $c(3,\delta) = 2^{-2^{C/\delta^C}}$ for some absolute constant $C$).

A more refined structure theorem was given in \cite{green-reg} (see also \cite{gt-arith}), which was termed an ``arithmetic regularity lemma'' in analogy with the Szemer\'edi regularity lemma which we discuss in the next section.  That theorem has similar hypotheses to Theorem \ref{stru}, but instead of constructing a single progression on $P$ on which one has pseudorandomness, one partitions $[1,N/3]$ into \emph{many} long progressions\footnote{Actually, for technical reasons it is more efficient to replace the notion of an arithmetic progression by a slightly different object known as a \emph{Bohr set}; see \cite{green-reg}, \cite{gt-arith} for details.}, where on most of which the function $f$ becomes linearly pseudorandom (after subtracting the mean).  A related structure theorem (with a more ``ergodic''
perspective) was also given in \cite{tao-survey}.  Here we give an alternate approach based on Fourier expansion and the pigeonhole principle.  Observe that for any $f: \ZBF/N\ZBF \to \CBF$ and any threshold $\lambda$ we have the Fourier decomposition $f = f_{U^\perp} + f_U$, where the ``structured'' component
$f_{U^\perp} := \sum_{\xi: |\hat f(\xi)| \geq \lambda} \hat f(\xi) e^{2\pi i x \xi/N}$ contains all the significant Fourier coefficients, and the ``pseudorandom'' component $f_U := \sum_{\xi: |\hat f(\xi)| \leq \lambda} \hat f(\xi) e^{2\pi i x \xi/N}$
contains all the small Fourier coefficients. Using Plancherel's theorem one can easily establish

\begin{theorem}[Weak structure theorem]\label{wst}  Let $f: \ZBF/N\ZBF \to \CBF$ be a function bounded in magnitude by $1$, and let $0 < \lambda < 1$.  Then we can split $f = f_{U^\perp} + f_U$, where $f_{U^\perp}$ is the linear combination of at most $O(1/\lambda^2)$ linear phase functions $x \mapsto e^{2\pi i x \xi/N}$, and $f_U$ is linearly pseudorandom in the sense that $\|f_U\|_{U^2(\ZBF/N\ZBF)} \leq \lambda$.
\end{theorem}

This theorem asserts that an arbitrary bounded function only has a bounded amount of significant linear Fourier-analytic structure; after removing this bounded amount of structure, the remainder is linearly pseudorandom.

This theorem, while simple to state and prove, has two weaknesses which make it unsuitable for such tasks as counting progressions of length three.  Firstly, even though $f$ is bounded by $1$, the components $f_{U^\perp}, f_U$ need not be.  Related to this,
if $f$ is non-negative, there is no reason why $f_{U^\perp}$ should be non-negative also.  Secondly, the pseudorandomness control
on $f_U$ is not very good when compared against the complexity of $f_{U^\perp}$ (i.e. the number of linear exponentials needed to describe $f_{U^\perp}$).  In practice, this means that any control one obtains on the structured component of $f$ will be dominated by the errors one has to concede from the pseudorandom component.  Fortunately, both of these defects can be repaired, the former by a Fej\'er summation argument, and the latter by a pigeonhole argument (which introduces a second error term $f_S$, which is small in $L^2$ norm).  More precisely, we have

\begin{theorem}[Strong structure theorem]\label{strong-struct}  Let $f: \ZBF/N\ZBF \to \RBF$ be a non-negative function bounded by $1$, and let $0 < \eps < 1$.  Let $F: \NBF \to \NBF$ be an arbitrary increasing function (e.g. $F(n) = 2^{2^n}$).  
Then there exists an integer $T = O_{F,\eps}(1)$ and a decomposition
$f = f_{U^\perp} + f_S + f_U$, where $f_{U^\perp}$ is the linear combination of at most $T$ linear phase functions, $f_U$ is linearly pseudorandom in the sense that $\|f_U\|_{U^2(\ZBF/N\ZBF)} = O(1/F(T))$, and $f_S$ is small in the sense that
$\| f_S\|_{L^2(\ZBF/N\ZBF)} := (\EBF_{n \in \ZBF/N\ZBF} |f_S(n)|^2)^{1/2} = O(\eps)$.  Furthermore, $f_{U^\perp}, f_U$ are bounded in magnitude by $1$.  Also, $f_{U^\perp}$ and $f_{U^\perp}+f_S$ are non-negative with the same mean as $f$.
\end{theorem}

\begin{proof} We use an argument from \cite{gk}.  We may take $\eps = 1/M$ for some large integer $M$.  
Let $N_1, N_2, \ldots, N_{M^2+2}$ be defined recursively by
$N_1 := M$ and $N_{m+1} := F( G(N_m) )^4$, where $G: \NBF \to \NBF$ is a function depending on $\eps$ 
that we shall choose later.  From Plancherel's theorem we have
$$ \sum_{\xi \in \ZBF/N\ZBF} |\hat f(\xi)|^2 \leq 1$$
and hence by the pigeonhole principle we can find $1 \leq m \leq M^2$ such that
$$ \sum_{ 1/N_{m+2} \leq |\hat f(\xi)| \leq 1/N_m} |\hat f(\xi)|^2 \leq 2/M^2 = O(\eps^2).$$
Now, for each $1 \leq m \leq M^2$, we define a Fej\'er-like kernel $K^{(m)}: \ZBF/N\ZBF \to \RBF^+$ which is 
non-negative, has mean one, has Fourier coefficients $1 + O(\eps)$ for all $\xi$ with $|\hat f(\xi)| \geq 1/N_m$,
and is a linear combination of at most $O_{N_m,\eps}(1)$ linear phase functions.  Such a function can be constructed
in a ``hard'' manner by means of Riesz products, or in a more ``soft'' manner by using the Weierstrass approximation theorem;
we omit the details.  If we then set
$$ f_{U^\perp} := f * K^{(m)}; \quad f_S := f * K^{(m+1)} - f * K^{(m)}; \quad f_U := f - f * K^{(m+1)},$$
with $T$ equal to the number of linear phase functions comprising $K^{(m)}$, then by repeated use of Plancherel's theorem 
one can verify all the required properties (if the function $G$ is chosen sufficiently fast growing, depending on $\eps$).
\end{proof}

Note that we have the freedom to set the growth function $F$ arbitrarily fast in the above proposition; this corresponds roughly
speaking to the fact that in the ergodic counterpart to this structure theorem (Proposition \ref{struct-3}) the pseudorandom error $f_U$ has asymptotically \emph{vanishing} Gowers $U^2$ norm.  One can view $f_{U^\perp}$ as a ``coarse'' Fourier approximation to $f$, and $f_{U^\perp}+f_S$ as a ``fine'' Fourier approximation to $f$; this perspective links this proposition with the graph regularity lemmas that we discuss in the next section.

Theorem \ref{strong-struct} can be used to deduce the structure theorems
in \cite{green-reg}, \cite{tao-survey}, \cite{gt-arith}, while a closely related result was also established in \cite{bourgain-firstroth}.  It can also be used to directly derive the $k=3$ case of Theorem \ref{sz-quant}, as follows.  Let $f$ be as in that proposition, and let $\eps := \delta^3/100$.  We apply Theorem \ref{strong-struct} to decompose $f = f_{U^\perp} + f_S + f_U$.  Because $f_{U^\perp}$ has only $T$ Fourier exponentials, it is easy to see that $f_{U^\perp}$ is almost periodic, in the sense that $\| T^n f_{U^\perp} - f_{U^\perp} \|_{L^2(\ZBF/N\ZBF)} \leq \eps$ for at least $c(\eps,T) N$ values of $n \in \ZBF/N\ZBF$, for some $c(\eps,T) > 0$.  For such values of $n$, one can easily verify that
$$ \EBF_{x \in \ZBF/N\ZBF} f_{U^\perp}(x) T^n f_{U^\perp}(x) T^{2n} f_{U^\perp}(x) \geq \EBF_{x \in \ZBF/N\ZBF} f_{U^\perp}^3 - 3\eps
\geq (\EBF_{x \in \ZBF/N\ZBF} f_{U^\perp})^3 - 3\eps \geq \delta^3/2.$$
Because $f_S$ is small, we can also deduce that
$$ \EBF_{x \in \ZBF/N\ZBF} (f_{U^\perp}+f_S)(x) T^n (f_{U^\perp}+f_S)(x) T^{2n} (f_{U^\perp}+f_S)(x) \geq \delta^3/4$$
for these values of $n$.  Averaging in $n$ (and taking advantage of the non-negativity of $f_{U^\perp} + f_S$) we conclude that
$$ \EBF_{x,n \in \ZBF/N\ZBF} (f_{U^\perp}+f_S)(x) T^n (f_{U^\perp}+f_S)(x) T^{2n} (f_{U^\perp}+f_S)(x) \geq\delta^3 c(\eps,T)/4.$$
Adding in the pseudorandom error $f_U$ using the generalized von Neumann theorem \eqref{gvn}, we conclude that
$$ \EBF_{x,n \in \ZBF/N\ZBF} f(x) T^n f(x) T^{2n} f(x) \geq  \delta^3 c(\eps,T)/4 - O( 1 / F(T) ).$$
If we choose $F$ to be sufficiently rapidly growing depending on $\delta$ and $\eps$, we can absorb the error term in the main term
and conclude that
$$ \EBF_{x,n \in \ZBF/N\ZBF} f(x) T^n f(x) T^{2n} f(x) \geq  \delta^3 c(\eps,T)/8.$$
Since $T = O_{F,\eps}(1) = O_{\delta}(1)$, we obtain the $k=3$ case of Theorem \ref{sz-quant} as desired.

Roth's original Fourier-analytic argument was published in 1953.  But the extension of this Fourier argument to the $k > 3$ case was not achieved until the work of Gowers \cite{gowers-4}, \cite{gowers} in 1998.  For simplicity we once again restrict attention to the $k=4$ case, where the theory is more complete.  Our objective is to show
\begin{equation}\label{quadruple-fourier}
\EBF_{x,r \in \ZBF/N\ZBF} f(x) T^r f(x) T^{2r} f(x) T^{3r} f(x) \geq c(4,\delta) > 0
\end{equation}
whenever $f$ is non-negative, bounded by $1$, and has mean at least $\delta$.  There are some significant differences between this case and the $k=3$ case \eqref{triple-fourier}.  Firstly, linear pseudorandomness is not enough to guarantee that a contribution to \eqref{quadruple-fourier} is negligible: for instance, if $f(x) := e^{2\pi i \xi x^2 / N}$, then
$$ \EBF_{x,r \in \ZBF/N\ZBF} f(x) T^r \overline{f}^3(x) T^{2r} f^3(x) T^{3r} \overline{f}(x) = 1$$
despite $f$ being very linearly pseudorandom (the $U^2$ norm of $f$ is $N^{-1/4}$); compare this example with \eqref{erg-ex}.  
One must now utilize some sort of ``quadratic Fourier analysis'' in order to capture the correct concept of pseudorandomness and structure.  Secondly, the Fourier-analytic arguments must now be supplemented by some results from arithmetic combinatorics (notably the Balog-Szemer\'edi theorem, and results related to Freiman's inverse sumset theorem) in order to obtain a usable 
notion of quadratic structure.  Finally, as in the ergodic case, one cannot rely purely on quadratic phase functions such as $e^{2\pi i (\xi x^2 + \eta x)/N}$ to generate all the relevant structured objects, and must also consider generalized quadratic objects such as locally quadratic phase functions, $2$-step nilsequences (see below), or bracket quadratic phases such as $e^{2\pi i \lfloor \sqrt{2} n \rfloor \sqrt{3} n}$.  

Let us now briefly sketch how the theory works in the $k=4$ case.  The correct notion of pseudorandomness is now given by the \emph{Gowers $U^3$ uniformity norm}, defined by
$$ \| f \|_{U^3(\ZBF/N\ZBF)}^8 := \EBF_{n \in \ZBF/N\ZBF} \| T^n f \overline{f} \|_{U^2(\ZBF/N\ZBF)}^4.$$
This norm measures the extent to which $f$ behaves quadratically; for instance, if $f = e^{2\pi i P(x)/N}$ for some polynomial $P$ of degree $k$ in the finite field $\ZBF/N\ZBF$, then one can verify that $\|f\|_{U^3(\ZBF/N\ZBF)}=1$ if $P$ has degree at most $2$, but (using the Weil estimates) we have $\|f\|_{U^3(\ZBF/N\ZBF)} = O_k(N^{-1/16})$ if $P$ has degree $k > 2$.  Repeated application of Cauchy-Schwarz then yields the generalized von Neumann theorem
\begin{equation}\label{gvn-4}
|\EBF_{x,r \in \ZBF/N\ZBF} f_0(x) T^r f_1(x) T^{2r} f_2(x) T^{3r} f_3(x)| \leq \min_{0 \leq j \leq 3}  \|f_j\|_{U^3(\ZBF/N\ZBF)}
\end{equation}
whenever $f_0,f_1,f_2,f_3$ are bounded in magnitude by $1$.
The next step is to establish a dichotomy between quadratic structure and quadratic pseudorandomness in the spirit of
Lemma \ref{lemcor-planch}.  In the original work of Gowers \cite{gowers-4}, it was shown that a function which was not quadratically pseudorandom had local correlation with quadratic phases on medium-length arithmetic progressions.  This result (when combined with the density increment argument of Roth) was already enough to prove \eqref{quadruple-fourier} with a reasonable
bound on $c(4,\delta)$ (basically of the form $1/\exp(\exp(\delta^{-C}))$); see \cite{gowers-4}, \cite{gowers}.  Building upon this work, a stronger dichotomy, similar in spirit to Lemma \ref{lemcor-4}, was established in \cite{gt-qm}.  Here, a number of essentially equivalent formulations of quadratic structure were established, but the easiest to state (and the one which generalizes most easily to higher $k$) is that of a \emph{(basic) $2$-step nilsequence}, which can be viewed as a notion of ``quadratic almost periodicity'' for sequences.  More precisely, a $2$-step nilsequence a sequence of the form $n \mapsto F(T^n x \Gamma)$, where $F$ is a Lipschitz function on a $2$-step nilmanifold $G/\Gamma$, $x\Gamma$ is a point in this nilmanifold, and $T$ is a shift operator $T: x\Gamma \mapsto \alpha x \Gamma$ for some fixed group element $\alpha \in G$.  We remark that quadratic phase sequences such as $n \mapsto e^{2\pi i \alpha n^2}$ are examples of $2$-step nilsequences, and generalized quadratics
such as $n \mapsto e^{2\pi i \lfloor \sqrt{2} n \rfloor \sqrt{3} n}$ can also be written (outside of sets of arbitrarily small density) as $2$-step nilsequences.

\begin{lemma}[Dichotomy between randomness and structure]\label{lemcor-planch4}\cite{gt-qm} Suppose that $f: \ZBF/N\ZBF \to \CBF$ is bounded in magnitude by $1$ with $\|f\|_{U^3(\ZBF/N\ZBF)} \geq \eta$ for some $0 < \eta < 1$.  Then there exists a $2$-step nilsequence $n \mapsto F(T^n x \Gamma)$, where $G/\Gamma$ is a nilmanifold of dimension $O_\eta(1)$, and $F$ is a bounded Lipschitz function $G/\Gamma$ with Lipschitz constant $O_\eta(1)$, such that $|\EBF_{1 \leq x \leq N} f(x) \overline{F(T^n x \Gamma)}| \geq c(\eta)$ for some $c(\eta) > 1$.  (We identify the integers from $1$ to $N$ with $\ZBF/N\ZBF$ in the usual manner.)
\end{lemma}

In fact the nilmanifold $G/\Gamma$ constructed in \cite{gt-qm} is of a very explicit form, being the direct sum of at most $O_\eta(1)$ circles (which are one-dimensional), skew shifts (which are two-dimensional), and Heisenberg nilmanifolds (which are three-dimensional).  The dimension $O_\eta(1)$ is in fact known to be polynomial in $\eta$, but the best bounds for $c(\eta)$ are currently only exponential in nature. See \cite{gt-qm} for further details and discussion.

The proof of Lemma \ref{lemcor-planch4} is rather lengthy but can be summarized as follows.  If $f$ has large $U^3$ norm, then by definition $T^n f \overline{f}$ has large $U^2$ norm for many $n$.  Applying Lemma \ref{lemcor-planch}, this shows that for many $n$, $T^n f \overline{f}$ correlates with a linear phase function of some frequency $\xi(n)$ (which can be viewed as a kind of ``derivative'' of the phase of $f$ in the ``direction'' $n$).  Some manipulations involving the Cauchy-Schwarz inequality then show that $\xi(n)$ contains some additive structure (in that there are many quadruples $n_1,n_2,n_3,n_4$ with $n_1+n_2=n_3+n_4$ and
$\xi(n_1)+\xi(n_2) = \xi(n_3)+\xi(n_4)$).  Methods from additive combinatorics (notably the Balog-Szemer\'edi(-Gowers) theorem and
Freiman's theorem, see e.g. \cite{tao-vu}) are then used to ``linearize'' $\xi$, in the sense that $\xi(n)$ agrees with a (generalized) linear function of $n$ on a large (generalized) arithmetic progression.  One then ``integrates'' this fact to conclude that $f$ itself correlates with a certain ``anti-derivative'' of $\xi(n)$, which is a (generalized) quadratic function on this progression.  This in turn can be approximated by a $2$-step nilsequence.  For full details, see \cite{gt-qm}.

Thus, quadratic nilsequences are the only obstruction to a function being quadratically pseudorandom.  This can be iterated to obtain structural results.  The following ``weak'' structural theorem is already quite useful:

\begin{theorem}[Weak structure theorem]\label{quad-wst}\cite{gt-arith}  Let $f: \ZBF/N\ZBF \to \CBF$ be a function bounded in magnitude by $1$, and let $0 < \lambda < 1$.  Then we can split $f = f_{U^\perp} + f_U$, where $f_{U^\perp}$ is a $2$-step nilsequence given by a nilmanifold of dimension $O_\lambda(1)$ and by a bounded Lipschitz function $F$ with Lipschitz constant $O_\lambda(1)$, and $f_U$ is quadratically pseudorandom in the sense that $\|f_U\|_{U^3(\ZBF/N\ZBF)} \leq \lambda$.  Furthermore, $f_{U^\perp}$ is non-negative, bounded by $1$, and has the same mean as $f$.
\end{theorem}

This is an analogue of Theorem \ref{wst}, and asserts that any bounded function has only a bounded amount of quadratic structure, with the function becoming quadratically pseudorandom once this structure is subtracted.  It cannot be proven in quite the same way as in Theorem \ref{wst}, because we have no ``quadratic Fourier inversion formula'' that decomposes a function neatly into quadratic components (the problem being that there are so many quadratic objects that such a formula is necessarily overdetermined).  However, one can proceed by a finitary analogue of the ergodic theory approach, known as an ``energy increment argument''.  In the ergodic setting, one uses all the quadratic objects to create a $\sigma$-algebra ${\mathcal Z}_2$, and sets $f_{U^\perp}$ to be the conditional expectation of $f$ with respect to that $\sigma$-algebra.  In the finitary setting, it turns out to be too expensive to try to use \emph{all} the $2$-step nilsequences to create a $\sigma$-algebra.  However, by adopting a more adaptive approach, selecting only those $2$-step nilsequences which have some significant correlation with $f$ (or some component of $f$), one can obtain the above theorem as follows.

\begin{proof}(Sketch) We perform the following iteration procedure.  Initialize ${\mathcal Z}$ to be the trivial $\sigma$-algebra $\{ \emptyset, \ZBF/N\ZBF \}$.  If $f - \EBF( f | {\mathcal Z} )$ is already quadratically pseudorandom, then stop the iteration.  Otherwise, using Lemma \ref{lemcor-planch4} we know that $f - \EBF( f | {\mathcal Z} )$ correlates with some $2$-step nilsequence $g(n) = F(T^n x \Gamma)$.  We take the level sets of $g$ (suitably discretized) and add them to the $\sigma$-algebra ${\mathcal Z}$; the correlation of $f - \EBF( f | {\mathcal Z} )$ with $g$ ensures that the 
\emph{energy} $\| \EBF( f | {\mathcal Z} ) \|_{L^2}^2$ will increase significantly (by some amount $c(\eta) > 0$) when doing so; this is essentially Pythagoras' theorem.    Because $f$ is bounded by $1$, the energy cannot exceed $1$, and so the iteration will stop
after $O_\eta(1)$ steps.  When one does this, one obtains a splitting $f = \EBF(f|{\mathcal Z}) + (f - \EBF(f|{\mathcal Z}))$, where
$f - \EBF(f|{\mathcal Z})$ is quadratically pseudorandom, and $\EBF(f|{\mathcal Z})$ is the conditional expectation of $f$ with respect to a bounded number of $2$-step nilsequences.  By applications of Urysohn's lemma, the Weierstrass approximation theorem, and the fact that any polynomial combination of $2$-step nilsequences is again a $2$-step nilsequence, we can approximate $\EBF(f|{\mathcal Z})$ to arbitrary accuracy by a $2$-step nilsequence $f_{U^\perp}$ of bounded complexity; by being careful one can also ensure that $f_{U^\perp}$ remains non-negative and bounded by $1$.  Setting $f_U := f - f_{U^\perp}$ one obtains the claim.
\end{proof}

It is likely that quantitative versions of this structure theorem will improve the known bounds on Szemer\'edi's theorem 
in the $k=4$ case; see \cite{gt-r4-1}, \cite{gt-r4-2}, \cite{gt-r4-3}.  A closely related version of this 
argument was also essential in establishing
Theorem \ref{szemp-thm}, see Section \ref{primes-sec} below.

\section{Graph theory}

We now turn to the third major line of attack to Szemer\'edi's theorem, based on graph theory (and hypergraph theory), and which is perhaps the purest embodiment of the strategy of exploiting the dichotomy between randomness and structure.  For graphs, the relevant structure theorem is the \emph{Szemer\'edi regularity lemma}, which was developed in \cite{szemeredi} in the original proof of Szemer\'edi's theorem, and has since proven to have many further applications in graph theory and computer science; see \cite{komlos} for a survey.  More recently, the analogous regularity lemma for hypergraphs have been developed in \cite{gowers-reg}, \cite{nrs}, \cite{rs}, \cite{rodl}, \cite{rodl2}, \cite{tao:hyper}.  Roughly speaking, these very useful lemmas assert that any graph (binary relation) or hypergraph (higher order relation), no matter how complex, can be modelled effectively as a pseudorandom sub(hyper)graph of a finite complexity (hyper)graph.  Returning to the setting of the introduction, the graph regularity lemma would assert that there exists a colouring of the integers into finitely many colours such that relations such as $x-y \in A$ can be viewed approximately as pseudorandom relations, with the ``probability'' of the event $x-y \in A$ depending only on the colour of $x$ and $y$.

The strategy of the graph theory approach is to abstract away the arithmetic structure in Szemer\'edi's theorem, converting the problem to one of finding solutions to an abstract set of equations, which can be modeled by graphs or hypergraphs.  
As before, we first illustrate this with the simple case of the $k=3$ case of Szemer\'edi's theorem, which we will take in the form of Theorem \ref{sz-quant}.  For simplicity we specialize to the case when $f$ is the indicator function of a set $A$ (which thus has density at least $\delta$ in $\ZBF/N\ZBF$); it is easy to see (e.g. by probabilistic arguments) that this special case 
in fact implies the general case.  The key observation is that the problem of locating an arithmetic progression of length three can be recast as the problem of solving three constraints in three unknowns, where each constraint only involves two of the unknowns.  Specifically, if $x,y,z \in \ZBF/N\ZBF$ solve the system of constraints
\begin{equation}\label{system}
 \begin{array}{llll}
 & y &+ 2z &\in A \\
-x & & + z &\in A \\
-2x & - y & & \in A
\end{array}
\end{equation}
then $y+2z, -x+z, -2x-y$ is an arithmetic progression of length three in $A$.  Conversely, each such progression comes from exactly $N$ solutions to \eqref{system}.  Thus, it will suffice to show that there are at least $c(3,\delta)N^3$ solutions to \eqref{system}.  Note that we already can construct at least $\delta N^2$ ``trivial solutions'' to \eqref{system}, in which $y+2z = -x+z = -2x+y$ is an element of $A$.  Furthermore, these trivial solutions $(x,y,z)$ are ``edge-disjoint'' in the
sense that no two of these solutions share more than one value in common (i.e. if $(x,y,z)$ and $(x',y',z')$ are distinct trivial solutions then at most one of $x=x'$, $y=y'$, $z=z'$ are true).  It turns out that these trivial solutions automatically generate a large number of non-trivial solutions to \eqref{system} - without using any further arithmetic structure present in these constraints.  Indeed, the claim now follows from the following graph-theoretical statement.

\begin{lemma}[Triangle removal lemma]\label{tri-remove}\cite{rsz}  For every $0 < \delta < 1$ there exists $0 < \sigma < 1$ with the following property. Let $G = (V,E)$ be an (undirected) graph with $|V|=N$ vertices which contains fewer than $\sigma N^3$ triangles.  Then it is possible to remove $O(\delta N^2)$ edges from $G$ to create a graph $G'$ which contains no triangles whatsoever.  
\end{lemma}

To see how the triangle removal lemma implies the claim, consider a vertex set $V$ which consists of three copies $V_1, V_2, V_3$
of $\ZBF/N\ZBF$ (so $|V|=3N$), and consider the tripartite graph $G = (V,E)$ whose edges are of the form 
$$ E = \{ (y,z) \in V_2 \times V_3: y+2z \in A \} \cup \{ (x,z) \in V_1 \times V_3: -x+z \in A \} \cup
\{ (x,y) \in V_1 \times V_2: -2x - y \in A \}.$$
One can think of $G$ as a variant of the Cayley graph for $A$.  Observe that solutions to \eqref{system} are in one-to-one 
correspondence with triangles in $G$.  Furthermore, the $\delta N^2$ trivial solutions to \eqref{system} correspond to $\delta N^2$ edge-disjoint triangles in $G$.  Thus to delete all the triangles one needs to remove at least $\delta N^2$ edges.  Applying Lemma \ref{tri-remove} in the contrapositive (adjusting $N$, $\delta$, $\sigma$ by constants such as $3$ if necessary), we see that $G$ contains at least $\sigma N^3$ triangles for some $\sigma=\sigma(\delta) > 0$, and the claim follows.

The only known proof of the triangle removal lemma proceeds by a structure theorem for graphs known as the \emph{Szemer\'edi regularity lemma}.  In order to emphasize the similarities between this approach and the previously discussed approaches, we shall not use the standard formulation of this lemma, but instead use a more recent formulation from \cite{tao:revisited}, \cite{tao:hyper} (see also \cite{alon-shapira}, \cite{lovasz-szegedy}), which replaces graphs with functions, and then obtains a structure theorem decomposing such 
functions into a structured (finite complexity)
component, a small component, and a pseudorandom (regular) component.  More precisely, we work with functions $f: V \times V \to \RBF$; this can be thought of as a weighted, directed generalization of a graph on $V$ in which every edge $(x,y)$ is assigned a real-valued weight $f(x,y)$.  The first step is to define a notion of pseudorandomness.  For graphs, this concept is well understood.  There are many equivalent formulations of this concept (see  \cite{chung-graham-wilson}), but we shall adopt one particularly close to the analogous concepts in previous sections, by introducing the \emph{Gowers $\Box^2$ cube norm} as
$$ \|f\|_{\Box^2}^4 := \EBF_{x,y,x',y' \in V} f(x,y) f(x,y') f(x',y) f(x',y');$$
when $f$ is the incidence function of a graph, the right-hand side essentially counts the number of $4$-cycles in that graph.  Again, one can use the Cauchy-Schwarz inequality to establish that the $\Box^2$ norm is indeed a norm; alternatively, one can use spectral theory and observe that the $\Box^2$ norm is essentially the Schatten-von Neumann $p$-norm of $f$ with $p=4$.  We refer to $f$ as \emph{pseudorandom} if its $\Box^2$ norm is small.
By two applications of Cauchy-Schwarz we have the generalized von Neumann inequality
\begin{equation}\label{gvn-graph}
 | \EBF_{x,y,z \in V} f(x,y) g(y,z) h(z,x) | \leq \min( \|f\|_{\Box^2}, \|g\|_{\Box^2}, \|h\|_{\Box^2} )
 \end{equation}
whenever $f,g,h$ are bounded in magnitude by $1$ (note that this generalizes \eqref{triple-fourier}).  

The next step, as before, is to establish a dichotomy between pseudorandomness and structure.  The analogue of Lemma \ref{lemcor} or Lemma \ref{lemcor-planch} is

\begin{lemma}[Dichotomy between randomness and structure]\label{lemcor-graph} Suppose that $f: V \times V \to \RBF$ is bounded in magnitude by $1$ with $\|f\|_{\Box^2(\ZBF/N\ZBF)} \geq \eta$ for some $0 < \eta < 1$.  Then there exists sets $A, B \subset V$ such that
$|\EBF_{x,y \in V} f(x,y) 1_A(x) 1_B(y)| \geq \eta^4/4$.
Here $1_A(x)$ denotes the indicator function of $A$ (thus $1_A(x) = 1$ if $x \in A$ and $1_A(x) = 0$ otherwise).
\end{lemma}

\begin{proof} By the definition of $\Box^2$ and the pigeonhole principle, one can find $x',y'$ such that
$$ |\EBF_{x,y \in V} f(x,y) f(x,y') f(x',y) f(x',y')| \geq \eta^4.$$
By splitting $f(x,y')$ and $f(x',y)$ into positive and negative parts, we conclude that there exist non-negative functions $a(x), b(y)$ bounded by $1$ such that
$$ |\EBF_{x,y \in V} f(x,y) a(x) b(y)| \geq \eta^4/4.$$
Now letting $A$, $B$ be random subsets of $V$, with $x \in A$ and $y \in B$ holding with independent probabilities $a(x)$ and $b(y)$ respectively.  From linearity of expectation we see that the expected value of $\EBF_{x,y \in V} f(x,y) 1_A(x) 1_B(y)$
has magnitude at least $\eta^4/4$, and the claim follows.
\end{proof}

One can iterate this to obtain a weak version of the Szemer\'edi regularity lemma:

\begin{theorem}[Weak structure theorem]\label{wst-graph}\cite{frieze}  Let $f: V \times V \to \RBF$ be a non-negative function bounded by $1$, and let $\eps > 0$.  Then we can decompose $f = f_{U^\perp} + f_U$, where 
$f_{U^\perp} = \EBF(f | {\mathcal Z} \otimes {\mathcal Z})$, ${\mathcal Z}$ is a $\sigma$-algebra of $V$ generated by at most $2/\eps$ sets, and $\|f_U\|_{\Box^2} \leq \eps$.  
\end{theorem}

\begin{proof}  (Sketch) We perform the following ``energy increment argument'' iteration, as in Theorem \ref{quad-wst}.
Initialize ${\mathcal Z}$ to be the trivial $\sigma$-algebra $\{ \emptyset, V \}$ on $V$, thus the tensor product 
${\mathcal Z} \otimes {\mathcal Z}$ is the trivial $\sigma$-algebra on $V \times V$.  If $f - \EBF( f | {\mathcal Z} \otimes {\mathcal Z} )$ has a $\Box^2$ norm less than $\eps$, stop the iteration.  Otherwise, use Lemma \ref{lemcor-graph} to find sets $A, B$ such that $1_A(x) 1_B(y)$ correlates
with $f - \EBF( f | {\mathcal Z} \otimes {\mathcal Z} )$.  One then adds $A$ and $B$ to the $\sigma$-algebra 
${\mathcal Z}$; the correlation of $f - \EBF( f | {\mathcal Z} )$ with $1_A(x) 1_B(y)$ ensures that the 
\emph{energy} $\| \EBF( f | {\mathcal Z} \otimes {\mathcal Z}) \|_{L^2}^2$ will increase significantly 
(by some amount $c(\eta) > 0$) when doing so; this is essentially Pythagoras' theorem.    Because $f$ is bounded by $1$, 
the energy cannot exceed $1$, and so the iteration will stop
after $O_\eta(1)$ steps.  When one does this, one obtains the desired splitting with
$f_{U} := f - \EBF(f|{\mathcal Z})$ and $f_{U^\perp} := \EBF(f|{\mathcal Z})$.
\end{proof}

As with Theorem \ref{wst}, the above theorem is too weak to be of much use, becase the control one has on the pseudorandomness of $f_U$ is fairly poor compared to the control on the complexity of $f_{U^\perp}$.  The following strong version of the regularity lemma is far more useful (compare with Theorem \ref{strong-struct}):

\begin{theorem}[Strong structure theorem]\label{sst}\cite{tao:revisited} Let $f: V \times V \to \RBF$ be a non-negative function bounded by $1$, and let $\eps > 0$.  Let $F: \NBF \to \NBF$ be an arbitrary increasing function (e.g. $F(n) = 2^{2^n}$).  
Then there exists an integer $T = O_{F,\eps}(1)$ and a decomposition
$f = f_{U^\perp} + f_S + f_U$, where $f_{U^\perp} = \EBF( f | {\mathcal Z} \otimes {\mathcal Z})$,  ${\mathcal Z}$ is generated by at most $T$ sets in $V$, $f_U$ is  pseudorandom in the sense that $\|f_U\|_{\Box^2} = O(1/F(T))$, and $f_S$ is small in the sense that
$\| f_S\|_{L^2(V \times V)} := (\EBF_{x,y \in V} |f_S(x,y)|^2)^{1/2} = O(\eps)$.  Furthermore, $f_{U^\perp}, f_U$ are bounded in magnitude by $1$.  Also, $f_{U^\perp}$ and $f_{U^\perp} + f_S$ are non-negative and bounded by $1$.
\end{theorem}

\begin{proof}  (Sketch) We repeat the energy increment argument from Theorem \ref{wst-graph}, but supplement it with an application of the pigeonhole principle.  Construct a sequence ${\mathcal Z}^{(0)} \subseteq {\mathcal Z}^{(1)} \subseteq \ldots$ of $\sigma$-algebras on $V$, with ${\mathcal Z}^(0)$ being the trivial algebra, and each ${\mathcal Z}^{(n+1)}$ formed by adding two sets $A, B$ to ${\mathcal Z}^{(n)}$ in such a way as to maximize the energy $E_{n+1} := \| \EBF( f | {\mathcal Z}^{(n+1)} \otimes {\mathcal Z}^{(n+1)} ) \|_{L^2}^2$.  From Pythagoras's theorem we see that the $E_n$ are increasing, but are also bounded between $0$ and $1$.  From the pigeonhole principle\footnote{Here we are exploiting a finitary version of the well-known fact that every bounded monotone sequence is convergent.  The finitary version is that if $E_n$ is an increasing sequence bounded above by $1$, $\eps > 0$, and $F: \NBF \to \NBF$, then there exists $n = O_{F,\eps}(1)$ such that $E_{n+F(n)} \leq E_n + \eps$.  This follows by defining a sequence $n_1, n_2, \ldots$ recursively by $n_1 := 1$ and $n_{i+1} := n_i + F(n_i)$ and observing from the pigeonhole principle that $E_{n_{i+1}} \leq E_{n_i} + \eps$ for some $i = O(1/\eps)$.}, one can thus find a positive integer $n = O_{F,\eps}(1)$ such that
$E_{n + F(2n)^4+1} \leq E_n + \eps^2$.  A further application of the pigeonhole principle then allows us to find
$n \leq n' \leq n+F(2n)^4$ such that $E_{n'+1} \leq E_{n'} + 1/F(2^{4n})^4$.  We now set
$$ f_{U^\perp} := \EBF(f | {\mathcal Z}^{(n)} \otimes {\mathcal Z}^{(n)} ); \quad
f_S := \EBF(f | {\mathcal Z}^{(n')} \otimes {\mathcal Z}^{(n')} ) - \EBF(f | {\mathcal Z}^{(n)} \otimes {\mathcal Z}^{(n)} );
\quad f_U := f - \EBF(f | {\mathcal Z}^{(n')} \otimes {\mathcal Z}^{(n')} )$$
and ${\mathcal Z} := {\mathcal Z}^{(n)}$.  Since $E_{n'} \leq E_n + \eps^2$, we see from Pythagoras' theorem that $f_S$ has an $L^2$ norm of $O(\eps)$.  Finally, since $E_{n'+1} \leq E_{n'} + 1/F(2n)^4$, the arguments in Theorem \ref{wst-graph} give
$\|f_U\|_{\Box^2} = O(1/F(2n))$.  Setting $T := 2n$ we obtain the claim.
\end{proof}

We remark that one could also prove Theorem \ref{sst} by a technique more similar to that used to prove Theorem \ref{strong-struct} by viewing $f$ as a matrix and using its singular value decomposition (or eigenvalue decomposition, if $f$ is symmetric) as a substitute for the Fourier inversion formula.  We omit the details.  One can view $f_{U^\perp}$ as a ``coarse'' approximation to $f$, as it is measurable with respect to a fairly low-complexity $\sigma$-algebra, and $f_{U^\perp} + f_S = \EBF( f | {\mathcal Z}^{(n')} \otimes {\mathcal Z}^{(n')} )$ as a ``fine'' approximation to $f$, which is considerably more complex but is also a far better approximation to $f$, in fact the accuracy of the fine approximation exceeds the complexity of the coarse approximation by any specified growth function $F$.  Also the difference between the coarse and fine approximations is controlled by an arbitrarily smal constant $\eps$.

Theorem \ref{sst} already easily implies the Szemer\'edi regularity lemma in its traditional formulation; see \cite{tao:revisited}.
It also implies Lemma \ref{tri-remove}, similar to how Theorem \ref{strong-struct} implies the $k=3$ version of Szemer\'edi's theorem.  We sketch the proof as follows.  Set $f$ to be the indicator function of $G$, thus
\begin{equation}\label{fff}
\EBF_{x,y,z \in V} f(x,y) f(y,z) f(z,x) \leq \sigma.
\end{equation}
Apply Theorem \ref{sst} to
obtain a decomposition $f = f_{U^\perp} + f_S + f_U$, where $F$ and $\eps$ are to be chosen later.  The $\sigma$-algebra ${\mathcal Z}$ is generated by at most $T$ sets, and thus has at most $2^T$ atoms.  We now use this decomposition to remove some ``irregular'' components of $G$.  First we remove from $G$ all edges with at least one vertex lying in an atom which is ``small'' in the sense that its cardinality is less than $\delta N / 2^T$; this costs us at most $O(\delta N^2)$ edges.  We also remove from $G$ all edges connecting a pair of atoms $A,B$ on which $f_S$ is ``large'' in the sense that $\EBF_{x \in A, y \in B} |f_S(x,y)|^2 \geq \eps^2/\delta$; this also costs us at most $O(\delta N^2)$ edges.  Finally, we remove from $G$ all edges connecting a pair of atoms $A,B$ on which $f_{U^\perp}$ is smaller than $\delta$ (or equivalently, $\EBF_{x \in A, y\in B} f(x,y) \leq \delta$); this also
costs us $O(\delta N^2)$ edges.  After all these removals, the only pairs of atoms $A,B$ which still contribute to the reduced graph
$G'$ are those which are large (so that $|A|, |B| \geq \delta N/ 2^T$), on which $f_{U^\perp}$ is larger than $\delta$, 
and on which $|f_S|^2$ has mean less than $\eps^2/\delta$.  Let us call such pairs $(A,B)$ ``good''.  

Now suppose that this reduced graph $G'$ still contains at least one triangle.  Then there must be three atoms $A,B,C$ such that
all three pairs $(A,B)$, $(B,C)$, $(C,A)$ are good.  In particular from the largeness of $f_{U^\perp}$ we have
$$ \EBF_{x \in A, y\in B, z \in C} f_{U^\perp}(x,y) f_{U^\perp}(y,z) f_{U^\perp}(z,x) \geq \delta^3$$
and then by the smallness of $f_S$ we have
$$ \EBF_{x \in A, y\in B, z \in C} (f_{U^\perp}+f_S)(x,y) (f_{U^\perp}+f_S)(y,z) (f_{U^\perp}+f_S)(z,x) \geq \delta^3
- O( \eps^2 / \delta )$$
and thus by the largeness of $A,B,C$ and the non-negativity of $f_{U^\perp}+f_S$
$$ \EBF_{x,y,z \in V} (f_{U^\perp}+f_S)(x,y) (f_{U^\perp}+f_S)(y,z) (f_{U^\perp}+f_S)(z,x) \geq [\delta^3
- O( \eps^2 / \delta )] \delta^3 / 2^{3T}.$$
Now by by the generalized von Neumann theorem \eqref{gvn-graph} and the pseudorandomness of $f_U$ we have
$$ \EBF_{x,y,z \in V} f(x,y) f(y,z) f(z,x) \geq [\delta^3 - O( \eps^2 / \delta )] \delta^3 / 2^{3T} - O( 1 / F(T) ).$$
If we choose $\eps$ to be a small multiple of $\delta^2$, and $F(T)$ to be a large multiple of $2^{3T}/\delta^6$, we thus have
$$ \EBF_{x,y,z \in V} f(x,y) f(y,z) f(z,x) \geq \frac{1}{2} \delta^6 / 2^{3T} \geq c(\delta)$$
for some $c(\delta) > 0$ (since $T = O_{F,\eps}(1) = O_\delta(1)$).  This will contradict \eqref{fff} if $\sigma$ is sufficiently small.  Thus $G'$ does not contain any triangles, and we are done.

As in the other two approaches, the above arguments extend (with some additional difficulties) to higher values of $k$.  Again we restrict attention to the $k=4$ case for simplicity.  To locate a progression of length four in a set $A \subset \ZBF/N\ZBF$ is now equivalent to solving the system of constraints
\begin{equation}\label{system-4}
 \begin{array}{lllll}
 & y &+ 2z &+3w &\in A \\
-x & & + z & + 2w&\in A \\
-2x & - y & & + w& \in A \\
-3x & -2y & -z & & \in A.
\end{array}
\end{equation}
This in turn follows from a hypergraph analogue of the triangle removal lemma. Define a \emph{$3$-uniform hypergraph} to be
a pair $H = (V,E)$ where $V$ is a finite set of vertices and $E$ is a finite set of unordered triplets $(x,y,z)$ in $V$, which we refer to as the \emph{edges} of $H$.  Define a \emph{tetrahedron} in $H$ to be a quadruple $(x,y,z,w)$ of vertices such that all four triplets $(x,y,z), (y,z,w), (z,w,x), (w,x,y)$ are edges of $H$.

\begin{lemma}[Tetrahedron removal lemma]\label{tetra-remove}\cite{frankl}  For every $0 < \delta < 1$ there exists $0 < \sigma < 1$ with the following property. Let $H = (V,E)$ be a $3$-uniform hypergraph graph with $|V|=N$ vertices which contains fewer than $\sigma N^4$ tetrahedra.  Then it is possible to remove $O(\delta N^3)$ edges from $H$ to create a hypergraph $H'$ which contains no tetrahedra whatsoever.  
\end{lemma}

Letting $f$ be the indicator function of $H$, we now have a situation where
$$ \EBF_{x,y,z,w \in V} f(x,y,z) f(y,z,w) f(z,w,x) f(w,x,y) \leq \sigma$$
and we need to remove some small components from $f$ so that this average now vanishes completely.  Again, the key step here is
to obtain a structure theorem that decomposes $f$ into structured parts, small errors, and pseudorandom errors.  The notion of pseudorandomness is now captured by the Gowers $\Box^3$ cube norm, defined by
$$ \| f \|_{\Box^3}^8 := \EBF_{x,y,z,x',y',z' \in V} f(x,y,z) f(x,y,z') f(x,y',z) f(x,y',z') f(x',y,z) f(x',y,z') f(x',y',z) f(x',y',z')$$
which in the case when $f$ is the indicator function of a hypergraph $H$, is essentially counting the number of octahedra present in $H$.  One can obtain a strong structure theorem analogous to Theorem \ref{sst}, but with one significant difference.
In Theorem \ref{sst}, the structured component $f_{U^\perp}(x,y)$ can be broken up into a small number of components which are of the form 
$1_A(x) 1_B(y)$.  In the $3$-uniform hypergraph analogue of Theorem \ref{sst}, the structured component $f_{U^\perp}(x,y,z)$ will be broken up into a small number of components of the form $1_A(x,y) 1_B(y,z) 1_C(z,x)$.  It turns out that in order to conclude the proof of Lemma \ref{tetra-remove}, this structural decomposition is not sufficient by itself; one must also turn to the functions $1_A(x,y)$, $1_B(y,z)$, $1_C(z,x)$ generated by this structure theorem and decompose them further, essentially by 
invoking Theorem \ref{sst}.  This leads to some technical complications in the argument, although this approach to Szemer\'edi's theorem is still the most elementary and self-contained.  See \cite{gowers-reg}, \cite{nrs}, \cite{rs}, \cite{rodl}, \cite{rodl2}, \cite{tao:hyper} for details.

\section{The primes}\label{primes-sec}

Having surveyed the three major approaches to Szemer\'edi's theorem, we now turn to the question of counting progressions in the primes (or in dense subsets of the primes).  The major new difficulty here, of course, is that the primes have asymptotically zero density rather than positive density, and even the most recent quantitative bounds on Szemer\'edi's theorem (see the discussion after Theorem \ref{sz-quant}) are not strong enough by themselves to overcome the ``thinness'' of the primes.   However, it turns out that the primes (and functions supported on the primes) are still within the range of applicability of structure theorems.  For instance, to oversimplify dramatically, the structure theorem in \cite{gt-primes} essentially\footnote{This is a gross oversimplification.  The precise statement is that after eliminating obvious irregularities in the primes caused by small residue classes, and excluding a small and technical exceptional set, a normalized counting function on the primes can be decomposed as a bounded function (which is thus spread out over a set of positive density), plus a pseudorandom error.  Ignoring the initial elimination of obvious irregularities and the exceptional set, and pretending the bounded function was the indicator function of a positive density set $A$, one recovers the interpretation of the primes as a sparse pseudorandom subset of $A$.} represents the primes (or any dense subset of the primes) as a (sparse) pseudorandom subset of a set of positive density.  Since sets of positive density already contain many progressions thanks to Szemer\'edi's theorem, it turns out that enough of these progressions 
survive when passing to a pseudorandom subset that one can conclude Theorem \ref{sz-quant}.

Interestingly, Theorem \ref{szemp-thm} can be tackled by (quantitative) ergodic methods, by Fourier-analytic methods, and
by graph-theoretic methods, with the three approaches leading to slightly different results.  For instance,
the establishment of infinitely many progressions of length three in the primes by van der Corput \cite{van-der-corput}
was Fourier-analytic, as was the corresponding statement for dense subsets of the primes (i.e. the $k=3$ case
of Theorem \ref{szemp-thm}), proven 76 years later by Green \cite{green}.  The argument in \cite{gt-primes} which proves Theorem \ref{szemp-thm} in full combines ideas from all three approaches, but is closest in spirit to the ergodic approach, albeit set in the finitary context of a cyclic group $\ZBF/N\ZBF$ rather than on an infinitary measure space.  The argument in \cite{tao:gaussian}, which shows that the Gaussian primes (or any dense subset thereof) contains infinitely many constellations
of any prescribed shape, and can be viewed as a two-dimensional analogue of Theorem \ref{szemp-thm}, was proven via the (hyper)graph-theoretical approach.  Finally, a more recent argument in \cite{gt-mobius}, \cite{gt-primeasymptotic}, in which precise asymptotics for the number of progressions of length four in the primes are obtained, as well as a ``quadratic pseudorandomness'' estimate on a renormalized counting function for the primes, proceeds by returning back to the original Fourier-analytic approach, but now using quadratic Fourier-analytic tools (Lemma \ref{lemcor-planch4} and Theorem \ref{quad-wst}) rather than linear ones.

As mentioned in the introduction, these results are discussed in other surveys \cite{kra-survey}, \cite{green-survey}, \cite{tao-coates}, \cite{tao-survey}, \cite{host-survey}, and we will only sketch some highlights here.  In all the results, the
strategy is to try to isolate the ``structured'' component of the primes from the ``pseudorandom'' component.  There is some obvious structure present in the primes; for instance, they are almost all odd, they are almost all coprime to three, and so 
forth.  This obvious structure can be normalized away fairly easily. For instance, to remove the bias the primes have towards being odd, one can replace the primes $P = \{2,3,5,\ldots\}$ with the renormalized set $P_{2,1} := \{ n: 2n+1 \hbox{ prime}\} = \{1,2,3,5,\ldots\}$.  Each arithmetic progression in $P_{2,1}$ clearly induces a corresponding progression in $P$, but the set $P_{2,1}$ has no bias modulo $2$.  More generally, to reduce all the bias present in residue classes mod $p$ for all $p < w$ (where $w$ is a medium-sized parameter to be chosen later), one can work with a set $P_{W,b} := \{ n: Wn+b \hbox{ prime}\}$, where $W$ is the product of all the primes less than $w$ and $1 \leq b < W$ is a number coprime to $W$.  This ``$W$-trick'' allows for some technical simplifications.

Next, it is convenient not to work with the primes as a set, but rather as a renormalized counting function.  One convenient choice is the von Mangoldt function $\Lambda(n)$, defined as $\log p$ if $n$ is a power of a prime $p$ and $0$ otherwise.  Actually, because of the $W$-trick, it is better to consider a renormalized von Mangoldt function such as $\Lambda_{W,b}(n) := \frac{W}{\phi(W)} \Lambda(Wn+b)$, where $\phi(W)$ is the Euler totient function of $W$.  The prime number theorem in arithmetic progressions asserts that the asymptotic average value of $\Lambda_{W,b}(n)$ is equal to $1$.  To establish progressions of length $k$ in the primes, it suffices to obtain a nontrivial lower bound for the asymptotic value of the average
\begin{equation}\label{wbn}
\EBF_{1 \leq n,r \leq N} \Lambda_{W,b}(n) \Lambda_{W,b}(n+r) \ldots \Lambda_{W,b}(n+(k-1)r).
\end{equation}
In fact this quantity is conjectured to asymptotically equal $1$ as $W,N \to \infty$, with $W$ growing much slower than $N$ 
(a special case of the Hardy-Littlewood prime tuples conjecture); the intuition is that by removing all the bias present in the small residue classes, we have eliminated all the ``obvious'' structure in the primes, and the renormalized function $\Lambda_{W,b}$ should now fluctuate pseudorandomly around its mean value $1$.  However, this conjecture has only been verified in the cases $k=3,4$ (leading to an asymptotic count for the number of progressions of primes of length $k$ less than a large number $N$); for the cases $k>4$ we only have a lower bound of $c(k)$ for some small $c(k) > 0$.

Let us cheat slightly by pretending that $\Lambda_{W,b}$ is a function on the cyclic group $\ZBF/N\ZBF$ rather than on the integers $\ZBF$; there are some minor technical truncation issues that need to be addressed to pass from one to the other but we shall ignore them here.  In order to show that \eqref{wbn} is close to $1$, an obvious way to proceed would be to establish some kind of
pseudorandomness control on the deviation $\Lambda_{W,b}-1$ from the mean, and then some sort of generalized von Neumann theorem to show that this deviation is negligible.  Based on the experience with Szemer\'edi's theorem, one would expect linear pseudorandomness to be the correct notion for $k=3$, quadratic pseudorandomness for $k=4$, and so forth.  In the $k=3$ case it is indeed a standard computation (using Vinogradov's method, or a modern variant of that method such as the one based on Vaughan's identity) to show that $\Lambda_{W,b}-1$ is has small Fourier coefficients, which is a reasonable proxy for linear pseudorandomness; the point being that the $W$-trick has eliminated all the ``major arcs'' which would otherwise destroy the pseudorandomness.  It then remains to obtain
a generalized von Neumann theorem, similar to \eqref{gvn}.  In preceding sections, one was working with functions that were bounded
(and hence square integrable), and one could obtain these theorems easily from Plancherel's theorem.  In the current setting, the
$L^2$ estimates on $\Lambda_{W,b}$ are unfavourable, and what one needs instead is some sort of $l^p$ bound on the Fourier coefficients of $\Lambda_{W,b}$ for some $2 < p < 3$.  This can be done by a more careful application of Vinogradov's method, but can also be achieved using harmonic analysis methods arising from restriction theory; see \cite{green}, \cite{green-tao}.  The key new insight here is that while the Fourier coefficients of $\Lambda_{W,b}$ are difficult to understand directly, one can \emph{majorize} $\Lambda_{W,b}$ pointwise by (a constant multiple of) a much better behaved function $\nu$ of comparable size, whose Fourier coefficients are much easier to obtain bounds for (indeed $\nu$ is essentially linearly pseudorandom once one subtracts off its mean, which is essentially $1$).  This ``enveloping sieve'' $\nu$ is essentially the Selberg upper bound sieve, and can be viewed as a ``smoothed out'' version\footnote{What is essentially happening here is that we are viewing the primes not as a zero density subset of the integers, but as a positive density subset of a set of ``almost primes'' which can be controlled efficiently via sieve theory.} of $\Lambda_{W,b}$.  Restriction theory (related to the method of the large sieve) is then used to pass from Fourier control of $\nu$ to Fourier control of $\Lambda_{W,b}$. 

A similar idea was used in \cite{green}, \cite{green-tao} to establish the $k=3$ case of Theorem \ref{szemp-thm}; we sketch the argument from \cite{green-tao} here as follows.  The main objective is to establish a lower bound for expressions such as
\begin{equation}\label{lambda-3}
\EBF_{x,r \in \ZBF/N\ZBF} \Lambda_{W,b} 1_A(x) \Lambda_{W,b} 1_A(x+r) \Lambda_{W,b} 1_A(x+2r)
\end{equation}
for large sets $A$.
Restriction theory still allows us to obtain good $l^p$ upper bound for the Fourier coefficients of 
$\Lambda_{W,b} 1_A$.  This functions as a substitute for Plancherel's theorem (which is not favourable here), and one can now obtain structure theorems such as Theorem \ref{wst} (and with some more effort, Theorem \ref{strong-struct}).  This decomposes $\Lambda_{W,b} 1_A$ into some structured component $f_{U^\perp}$ and a linearly pseudorandom component $f_U$.  The generalized von Neumann theorem lets us dispose the contribution of $f_U$ to \eqref{lambda-3}, so let us focus on $f_{U^\perp}$.  One can try to use the complexity bound on $f_{U^\perp}$ (controlling the number of linear phases that comprise $f_{U^\perp}$) to get some lower bound here, but this would require developing a strong structure theorem analogous to Theorem \ref{strong-struct}.  It turns out that one can argue more cheaply, using a weaker structure theorem analogous to Theorem \ref{wst}.  The key observation is that because $\Lambda_{W,b} 1_A$ is dominated (up to a constant) by the enveloping sieve $\nu$, the structured component of $\Lambda_{W,b} 1_A$ (which is essentially a convolution of $\Lambda_{W,b} 1_A$ with a Fej\'er-like kernel) is pointwise dominated (up to a constant) by a corresponding structured component of $\nu$.  But since $\nu$ is linearly pseudorandom after subtracting off its mean, the structured component of $\nu$ turns out to essentially be just the mean of $\nu$, which is bounded.   We conclude that $f_{U^\perp}$ is bounded, at which point one can just apply Szemer\'edi's theorem (Theorem \ref{sz-quant}) directly to obtain a good lower bound on this contribution to \eqref{lambda-3}, and one can now conclude the $k=3$ case of Theorem \ref{szemp-thm}.

The proof of Theorem \ref{szemp-thm} for general $k$ in \cite{gt-primes} 
follows the same general strategy, but it is convenient to abandon the Fourier
framework (which becomes quite complicated for $k>3$) and instead take an approach which borrows ingredients from all three approaches, especially the ergodic theory approach.  
From the Fourier approach one borrows the Gowers uniformity norms $U^{k-2}(\ZBF/N\ZBF)$, which are a convenient way to define the appropriate notion of pseudorandomness for counting progressions of length $k$.  
One still needs an enveloping sieve $\nu$, but instead of using a Selberg-type sieve that enjoys good Fourier coefficient control, it turns out to be more convenient to use an enveloping sieve\footnote{A related enveloping sieve was also used in the recent establishment of narrow gaps in the primes \cite{pintz}.}
 of Goldston and Y{\i}ld{\i}r{\i}m \cite{goldston-yildirim-old1}, \cite{goldston-yildirim-old2}, \cite{goldston-yildirim} which has good control on $k$-point correlations (indeed, it behaves pseudorandomly after subtracting off its mean, which is essentially $1$).

The next step is a generalized von Neumann theorem to show that the contribution of pseudorandom functions are negligible.  The fact that the functions involved are no longer bounded by $1$, but are instead dominated by $\nu$, makes this theorem somewhat trickier to establish, however it can still be achieved by a number of applications of the Cauchy-Schwarz and taking advantage of the pseudorandomness properties of $\nu-1$.  This type of argument is inspired by certain ``sparse counting lemmas'' arising from the hypergraph approach, particuarly from \cite{gowers-reg}.

The main step, as in previous sections, is a structure theorem which decomposes $\Lambda_{W,b}$ (or $\Lambda_{W,b} 1_A$) into
a structured component and a pseudorandom component.  In principle one could use higher order Fourier analysis (or the
precise characteristic factors achieved in \cite{host-kra2}, \cite{ziegler} to obtain this decomposition, but this looks rather
difficult technically, though progress has been made in the $k=4$ case.  Fortunately, there is a ``softer'' approach in which one defines structure purely by duality; to oversimplify substantially, one defines a function to be structured if it is approximately orthogonal to all pseudorandom functions. One can then obtain a soft structural theorem in which the structural component is essentially a conditional expectation of the original function to a certain $\sigma$-algebra generated by certain special structured functions which are called ``dual functions'' in \cite{gt-primes}.  This $\sigma$-algebra (the finitary analogue of a characteristic factor) is not too tractable to work with, but somewhat miraculously, one can utilize the pseudorandomness properties of $\nu$ and a large number of applications of the Cauchy-Schwarz inequality to show that the conditional expectation of $\nu$ with respect to this $\sigma$-algebra remains bounded (outside of a small exceptional set, which turns out to have a negligible impact).  Since $\Lambda_{W,b} 1_A$ is pointwise dominated by a 
constant multiple of $\nu$, the structured component of $\Lambda_{W,b} 1_A$ is similarly bounded and can thus be controlled using Szemer\'edi's theorem.  Combining this with the generalized von Neumann theorem to handle the pseudorandom component, one obtains Theorem \ref{szemp-thm}.  The result for the Gaussian prime constellations is similar, but uses the Gowers cube norms $\Box^{k-2}$ instead of the uniformity norms, and replaces Szemer\'edi's theorem by a hypergraph removal lemma similar to Lemma \ref{tri-remove} and Lemma \ref{tetra-remove}; see \cite{tao:hyper}, \cite{tao:gaussian}.

The arguments used to prove Theorem \ref{szemp-thm} give a lower bound for the expression \eqref{wbn}, but do not compute its asymptotic value (which should be $1$).  As mentioned earlier, for $k=3$ this can be achieved by the circle method.  More recently, the $k=4$ case has been carried out in \cite{gt-mobius}, \cite{gt-primeasymptotic}; the same method in fact allows one to asymptotically count the number of solutions to any two linear homogeneous equations in four prime unknowns.  The key point is to show that $\Lambda_{W,b}-1$ is quadratically pseudorandom, as the generalized von Neumann theorem will then allow one to control \eqref{wbn} satisfactorily.  It turns out that a variant of Lemma \ref{lemcor-planch4} applies here, and reduces matters to showing that $\Lambda_{W,b}-1$ does not correlate significantly with any $2$-step nilsequences.  This task is attackable by Vinogradov's method, although it is rather lengthy and it turns out to be simpler to first replace $\Lambda_{W,b}-1$ with the closely related M\"obius function.

\end{document}